\documentclass[a4paper,11pt,oneside,reqno]{amsart}
\usepackage{amsmath,amsfonts,amsthm,amssymb,dsfont,tikz-cd}
\usepackage{geometry,enumitem,hyperref}
\title{Green points in the reals}
\author{Yilong Zhang}
\address{Mathematisches Institut, Universität Bonn, Endenicher Allee 60, D-53115 Bonn, Germany}
\email{s58yzhan@uni-bonn.de}
\date{\today}

\newenvironment{claim}[1]{\par\noindent\textbf{Claim.}\space#1}{}
\newenvironment{pfcl}[1]{\par\noindent\emph{Proof.}\space#1}{}

\theoremstyle{plain}
\newtheorem{theorem}{Theorem}[section]
\newtheorem{lemma}[theorem]{Lemma}
\newtheorem{proposition}[theorem]{Proposition}
\newtheorem{fact}[theorem]{Fact}
\newtheorem{corollary}[theorem]{Corollary}
\newtheorem{conjecture}[theorem]{Conjecture}
\newtheorem{maintheorem}{Theorem}

\theoremstyle{definition}
\newtheorem{definition}[theorem]{Definition}

\numberwithin{equation}{section}

\DeclareMathOperator{\cl}{cl}
\DeclareMathOperator{\Cl}{Cl}
\DeclareMathOperator{\rcl}{rcl}
\DeclareMathOperator{\dcl}{dcl}
\DeclareMathOperator{\trd}{trd}
\DeclareMathOperator{\ld}{ld}
\DeclareMathOperator{\md}{md}
\DeclareMathOperator{\rcf}{RCF}
\DeclareMathOperator{\fin}{fin}
\DeclareMathOperator{\rich}{rich}
\DeclareMathOperator{\tp}{tp}
\DeclareMathOperator{\qftp}{qf-tp}
\DeclareMathOperator{\dprk}{dp-rk}
\DeclareMathOperator{\rk}{rk}
\DeclareMathOperator{\zar}{Zar}
\def\forkindep{\mathrel{\raise0.2ex\hbox{\ooalign{\hidewidth$\vert$\hidewidth\cr\raise-0.9ex\hbox{$\smile$}}}}}

\newcommand{\A}{\mathcal{A}}
\newcommand{\B}{\mathcal{B}}
\newcommand{\C}{\mathcal{C}}

\newcommand{\F}{\mathcal{F}}
\newcommand{\M}{\mathcal{M}}
\newcommand{\I}{\mathcal{I}}
\newcommand{\R}{\mathcal{R}}
\newcommand{\U}{\mathcal{U}}
\newcommand{\X}{\mathcal{X}}
\newcommand{\Y}{\mathcal{Y}}
\newcommand{\Z}{\mathcal{Z}}
\newcommand{\W}{\mathcal{W}}
\newcommand{\bN}{\mathbb{N}}
\newcommand{\bQ}{\mathbb{Q}}
\newcommand{\bR}{\mathbb{R}}
\newcommand{\bC}{\mathbb{C}}
\newcommand{\bZ}{\mathbb{Z}}
\newcommand{\ba}{\boldsymbol{a}}
\newcommand{\bb}{\boldsymbol{b}}
\newcommand{\bc}{\boldsymbol{c}}
\newcommand{\bx}{\boldsymbol{x}}
\newcommand{\by}{\boldsymbol{y}}
\newcommand{\bz}{\boldsymbol{z}}

\newcommand{\bw}{\boldsymbol{w}}
\newcommand{\bs}{\boldsymbol{s}}
\newcommand{\bt}{\boldsymbol{t}}
\newcommand{\bp}{\boldsymbol{p}}
\newcommand{\bq}{\boldsymbol{q}}
\newcommand{\ha}{\widehat{\ba}}
\newcommand{\hb}{\widehat{\bb}}

\newcommand{\Cf}{\C^{\fin}}
\newcommand{\ag}[1]{{\langle#1\rangle}}
\newcommand{\se}{\subseteq}

\begin{document}

\begin{abstract}
  We construct an expansion of a real closed field by a multiplicative subgroup adapting Poizat's theory of green points. Its theory is strongly dependent, and every open set definable in a model of this theory is semialgebraic. We prove that the real field with a dense family of logarithmic spirals, proposed by Zilber, satisfies our theory.
\end{abstract}

\maketitle

\section{Introduction}

In model theory, the Fraïssé limit is a common tool for constructing new structures from their substructures. In the search of new strongly minimal structures, Hrushovski \cite{Hr93} adjusts the Fraïssé limit by restricting the amalgamation to ``strong'' embeddings, constructing a counterexample to Zilber's dichotomy conjecture. This method is now called the Hrushovski construction.

Using the Hrushovski construction, Poizat \cite{Po01} obtains new stable expansions of algebraically closed fields, in particular, a rank $\omega\cdot 2$ expansion of an algebraically closed field of characteristic 0 by a multiplicative subgroup, whose theory is called the theory of green points. Baudisch, Hils, Martin-Pizarro, and Wagner \cite{Ba09} extend this work by constructing a rank 2 expansion, proving the existence of bad fields of characteristic 0.

Zilber \cite{Zi04} studies natural models of the theories constructed by Poizat. As a result, Caycedo and Zilber \cite{Ca14} find a model of Poizat's theory of green points, which is the expansion of the complex field by a dense family of logarithmic spirals. They suggest that similar methods could be applied to the study of expansions of the real field, which is implemented in this paper.

In order to implement Zilber's proposal in the situation of the real field, we need a real version of Poizat's theory of green points. To achieve that, we apply the Hrushovski construction to expansions of a real closed field by a divisible multiplicative subgroup obtaining a theory $T^{\rich}_\B$. We then show that it is the theory of the real field with a dense family of logarithmic spirals, the same predicate as Caycedo and Zilber's.

Let $\overline{\bR}:=(\bR,<,+,-,\cdot,0,1)$ be the real ordered field. As usual, we identify $\bR^2$ with $\bC$. Let $(\overline{\bR},G)$ be the expansion by the subset $G\se\bR^2$ given by:
  \[G=\exp(\epsilon\bR+Q)=\{e^{\epsilon t+s}\in\bR^2\mid t\in\bR,s\in Q\},\]
where $Q\se\bR$ is a non-trivial finite dimensional $\bQ$-vector space, and $\epsilon=1+\beta i$ for some $\beta\in\bR\setminus\{0\}$. In the following, SC$_K$ stands for the Schanuel property for a field $K$ (Conjecture \ref{con:SCK}). By Bays, Kirby, and Wilkie \cite{Ki10}, SC$_K$ holds for $K=\bQ(\beta i)$ for every exponentially transcendental $\beta$.

\begin{maintheorem}\label{thm:md}
  If SC$_K$ holds for $K=\bQ(\beta i)$, then $(\overline{\bR},G)$ is a model of $T^{\rich}_\B$.
\end{maintheorem}

See Section \ref{ch:md} for the proof of the above theorem. Consequently, $(\overline{\bR},G)$ inherits model-theoretic properties from $T^{\rich}_\B$. Comparing with the expansions of the real field by locally closed trajectories of linear vector fields, which by Miller \cite{Mi11} either possess d-minimality or define the set of integers, our model $(\overline{\bR},G)$ is an expansion by a dense and codense subset, hence not d-minimal. Nevertheless, we show that all its definable open sets are semialgebraic, and therefore it has an o-minimal open core. For neostability, we prove that it is strongly dependent.

More precisely, we apply the Hrushovski construction to the following setting:
\begin{itemize}
  \item a real closed field $\R=(R,<,+,-,\cdot,0,1)$,
  \item a divisible subgroup $G$ of the multiplicative group of the algebraically closed field $K:=R^2$,
  \item the class $\C$ of all such structures $(\R,G)$,
  \item the divisible hull operator $\cl:\mathcal{P}(K^\times)\to\mathcal{P}(K^\times)$ given by
    \[\cl(A):=\{b\in K^\times\mid b^n\in\ag{A}\text{ for some }n\in\bN\},\text{ and}\]
  \item the predimension function $\delta(\ba):=\trd(\ba)-\md(\cl(\ba)\cap G)$, for $\ba\se R^2$.
\end{itemize}

The transcendence degree of $\ba=\big((x_1,y_1),...,(x_n,y_n)\big)$ is counted as follows:
\[\trd(\ba):=\trd(x_1,y_1,...,x_n,y_n).\]
This is because we have the information on both real and imaginary axes in our language. The multiplicative degree $\md$ is the dimension function induced by $\cl$.

\emph{Strong embeddings} between structures in $\C$ are embeddings with non-negative predimension. Fix a structure $\B\in\C$ finitely rcl-generated (finitely generated as a real closed field). Let $\C_\B$ consist of all structures $\A\in\C$ such that $\B$ is strongly embedded into $\A$, and let $\Cf_\B$ consist of all finitely rcl-generated structures in $\C_\B$.

A structure $\U\in\C_\B$ is \emph{rich} if for every $\Y\in\Cf_\B$ and every strong substructure $\X$ of $\Y$, every strong embedding from $\X$ to $\U$ extends to a strong embedding from $\Y$ to $\U$. We show that $\C_\B$ has the amalgamation property, and hence contains a rich structure. Rich structures are back-and-forth equivalent.

In rich structures, the predimension function induces a pregeometry $\Cl$. When applied to expansions of strongly minimal structures, the Hrushovski construction produces stable structures, because the $\Cl$-independent type is unique. In the ordered case, we have the analogue that the $\Cl$-independent type for a fixed cut is unique. This property allows us to show that, instead of stability, every open set definable in a model of $T^{\rich}_\B$ is semialgebraic. In the terminology of Dolich, Miller, and Steinhorn \cite{Do10}, it implies that the theory of real closed fields is an open core of $T^{\rich}_\B$. We also prove that over a finite set $A$ the types inside $\Cl(A)$ have finite dp-rank. As a result, rich structures are strongly dependent.

To find the theory of rich structures, we need to first axiomatize the class $\C_\B$. The idea is that, if $\ba$ satisfies too many algebraic relations, then $\ba$ must belong to some coset of a proper algebraic subgroup of the torus. By Zilber's weak CIT (Fact \ref{fact:WCIT}), given a variety $W\se K^n$, there is a finite set $\mu(W)$ of proper algebraic subgroups such that, every atypical intersection of $W$ and a coset of an algebraic subgroup is contained in a coset of some subgroup in $\mu(W)$. Consider a variety $V\se R^{2n}$ over $\B$. Let $\widehat{V}$ be the Zariski closure of $V$ in $K^n$, and let $\widetilde{V}$ be the Zariski closure of $(V,V^c)$ in $K^{2n}$. Let $\phi_V$ be the sentence stating the following:
\begin{align*}
&\forall\,(x_1,y_1,...,x_n,y_n)\in V\cap G^n,\\
&\text{the fiber }\widetilde{V}(\bx+i\by,-)\text{ is generic }\longrightarrow\bigvee_{A\in\mu(\widehat{V})}\bx+i\by\in\beta A\text{ for some }\beta\in (G^\B)^n.
\end{align*}
For the precise definition of $\phi_V$, see Subsection \ref{sec:alo}. Let $T_\B$ contain $\phi_V$ for all $V$ of dimension $<n$. Then $T_\B$ axiomatizes $\C_\B$.

To capture existential closedness of rich structures, we use a set of axioms EC, saying that every rotund block intersects $G$, where rotund blocks are precisely the loci of strong extensions. Let $T_\B^{\rich}$ be the union of $T_\B$ and EC. We prove the following theorem in Subsection \ref{sec:rich}.

\begin{maintheorem}
  \label{thm:th}
  Rich structures in $\C_\B$ are exactly $\omega$-saturated models of $T_\B^{\rich}$.
\end{maintheorem}

The desired structure $(\overline{\bR},G)$ satisfies EC without further assumption. The key observation is that every rotund set has locally finite intersection with each logarithmic spiral. Thus, fiber dimension assures the intersection with a dense subset. In order to find a finitely rcl-generated substructure $\B$ such that $(\overline{\bR},G)\in\C_\B$, we need the Schanuel property (SC$_K$) for a lower bound of the predimension function.

The paper is structured as follows. In Section \ref{ch:pr}, we recall the Hrushovski construction and some results in algebra. Section \ref{ch:cons} is devoted to the construction of rich structures. In Section \ref{ch:rt}, we define rotundity and give its first-order characterization. In Section \ref{ch:rs}, we establish basic properties for rich structures. Section \ref{ch:th} is devoted to the axiomatization of rich structures. In Section \ref{ch:tm}, we verify near model completeness, o-minimal open core, and strong dependence of the theory. In Section \ref{ch:md}, we prove that $(\overline{\bR},G)$ is a model of the theory of rich structures.

\subsection*{Acknowledgements}

I would like to thank my supervisor, Prof. Philipp Hieronymi, for suggesting this topic and improving this paper. I am grateful to Prof. Martin Hils for his advice on the draft. I thank the Mathematical Institute of the University of Bonn for supporting me during the research of this work.

\subsection*{Notations}

Throughout this article, we work in first-order languages and one-sorted structures. We use curly letter $\X,\Y$ for structures and $X,Y$ for their underlying sets, and we use bold letters $\ba,\bb$ for finite tuples. We write $\ba\bb$ for a concatenation of two tuples, and we write $\tp$ for type, $\qftp$ for quantifier-free type. Let $(-)^c$ be the complex conjugation. Let $\ld_K,\md,\trd$ denote $K$-linear dimension for some field $K$, multiplicative degree, and transcendence degree respectively.

Let $S$ be a set. For a subset $X\se S^{n+m}$ and a tuple $\ba\in S^m$, we write $X(-,\ba)$ or $X_{\ba}$ for the \emph{fiber} of $\ba$, i.e.,
\[X(-,\ba)=X_{\ba}=\{\,\bx\in S^n\mid(\bx,\ba)\in X\,\}.\]

Let $\R\vDash\rcf$. Given a set $A\se R$, let $\rcl(A)$, or simply $\ag{A}$, denote the smallest real closed subfield of $R$ containing $A$. For a set $W\se R^n$, let $\dim^\R(W)$ denote the topological dimension of $W$. For a field $K$ and a variety $V\se K^n$, let $\dim^K(V)$ denote the Krull dimension of $V$. A point $\ba$ in $W$ (or in $V$) is \emph{generic} if its transcendence degree equals the dimension of $W$ (or of $V$).

Let $\cl$ a pregeometry. We denote $A\forkindep^{\cl}_C B$ if $A$ and $B$ are $\cl$-independent over $C$, i.e., for all finite $A_0\se A$,
\[\dim^{\cl}(A_0/BC)=\dim^{\cl}(A_0/C),\]
where $\dim^{\cl}$ is the dimension function induced by $\cl$.

\section{Preliminaries}
\label{ch:pr}

\subsection{Hrushovski construction}
\label{sec:hc}

In this subsection, we review the concepts related to Hrushovski's predimension construction method, and adapt them to our context, i.e., expansions of fields.

\subsubsection{Predimension}

Let $A$ be a set and let $\cl$ be a modular pregeometry on $A$. Consider a function $\delta:\mathcal{P}_{\fin}(A)\to\bZ$.

\begin{definition}
  A function $\delta$ is \emph{well-defined with respect to $\cl$} if for every finite tuple $\bx,\by$ with $\cl(\bx)=\cl(\by)$ we have $\delta(\bx)=\delta(\by)$.
\end{definition}

Any function well-defined with respect to $\cl$ can also be defined for finite dimensional $\cl$-closed subsets. Indeed, take $\delta(X)$ to be $\delta(\bx)$ whenever $X=\cl(\bx)$.

\begin{definition}
  \label{de:loc}
  The \emph{localization of $\delta$ at $\by$} is the function
  \[\delta_{\by}(\bx)=\delta(\bx/\by):=\delta(\bx\by)-\delta(\by).\]
\end{definition}

\begin{definition}
  A function $\delta$ well-defined with respect to $\cl$ is a \emph{predimension function} if for every finite dimensional $\cl$-closed $X$ and $Y$, the submodularity condition holds:
  \[\delta(\cl(X\cup Y))\leq\delta(X)+\delta(Y)-\delta(X\cap Y).\]
\end{definition}

Any localization of a predimension function is also a predimension function.

From now on, suppose $\delta$ is a predimension function with respect to $\cl$.

\begin{definition}
  \label{def:ssf}
  Let $X\se Y$ be $\cl$-closed subsets of $A$, where $X$ has finite $\cl$-dimension. We say $X$ is \emph{strong} in $Y$ if $\delta(\by/X)\geq 0$ for every $\by\se Y$.
\end{definition}

\begin{definition}
  \label{def:ss}
  Let $X\se Y$ be any $\cl$-closed subsets of $A$. We say $X$ is \emph{strong} in $Y$ if $X$ is the union of a directed set of finite dimensional $\cl$-closed subsets, all of which are strong in $Y$.
\end{definition}

Let us write $X\leq Y$ for $X$ being a strong subset of $Y$. We collect some basic facts about strong subsets in the following proposition. For the proof, we refer to \cite[Section 2.2.2]{Ca11}

\begin{fact}\label{fact:str}
  Let $X,Y,Z$ be $\cl$-closed subsets of $A$.
  \begin{enumerate}
    \item (Subset) If $X\se Y\se Z$ and $X\leq Z$, then $X\leq Y$.
    \item (Transitivity) If $X\leq Y$ and $Y\leq Z$, then $X\leq Z$.
    \item (Union) Let $(X_i)_{i\in I}$ be an ascending $\leq$-chain. Then $X_j\leq\bigcup_{i\in I}X_i$ for every $j\in I$.
    \item (Intersection) If $X\leq Z$ and $Y\leq Z$, then $X\cap Y\leq Z$.
  \end{enumerate}
\end{fact}

\begin{definition}
  Let $X\se A$. The \emph{hull} of $X$ in $A$, denoted by $[X]^A$, is the smallest $\cl$-closed strong subset of $A$ containing $X$. We omit the superscript $A$ when it is clear from the context.
\end{definition}

\begin{fact}
  \label{fact:min}
  Suppose $\delta:\mathcal{P}_{\fin}(A)\to\bZ$ is bounded below. Then $[X]$ exists for every $X\se A$. Furthermore, every $(x_1,...,x_n)\in A^k$ can be extended to a finite tuple $(x_1,...,x_n,y_1,...,y_m)$ such that $[\{x_1,...,x_n\}]=\cl(x_1,...,x_n,y_1,...,y_m)$.
\end{fact}

\begin{definition}
  The \emph{dimension function associated to $\delta$} is defined by
  \[d(\bx):=\min\{\delta(\by)\mid\bx\se\by\se A\}.\]
\end{definition}

If $\delta$ is bounded below, then $d(\bx)=\delta([\bx])$. Moreover, if $\delta$ satisfies:
\begin{itemize}
  \item $\delta(\bx)\geq 0$ for every $\bx\se A$,
  \item $\delta(\emptyset)=0$, and
  \item $\delta(\{x\})\leq 1$ for every $x\in A$,
\end{itemize}
then $d$ is a dimension function. We write $\Cl$ for the pregeometry given by $d$. It is easy to see that
\[\Cl(\bx)=\bigcup\{\by\se A\mid\delta(\by/[\bx])=0\}.\]

\subsubsection{Amalgamation}

Let $\C$ be a class of $L$-structures for some language $L$. Suppose we have a pregeometry and a predimension function defined uniformly on $\C$. More precisely, for each $\X\in\C$, we have a pregeometry $\cl|_{\X}$ and a predimension function $\delta|_{\X}$ with respect to $\cl|_{\X}$, both preserved under partial isomorphisms.

\begin{definition}
  Let $\X,\Y$ be structures in $\C$.

  If $\X$ is a substructure of $\Y$ and $X$ is strong in $Y$, then we say $\Y$ is a \emph{strong extension} of $\X$, denoted by $\X\leq\Y$.

  An embedding $\iota:\X\to\Y$ is a \emph{strong embedding} if $\iota(X)$ is strong in $Y$.
\end{definition}

Let $\C^{\fin}$ be the subclass of $\C$ consisting of all finitely generated structures.

\begin{definition}
  A structure $\U\in\C$ is \emph{rich} if for every $\X,\Y\in\C^{\fin}$ and strong embeddings $f:\X\to\Y,\,g:\X\to\U$, there is a strong embedding $h:\Y\to\U$ such that $h\circ f=g$.
  \begin{center}\begin{tikzcd}
    \X \arrow[rd, "g"] \arrow[d, "f"']\\
    \Y \arrow[r, "h", dashed] & \U
  \end{tikzcd}\end{center}
\end{definition}

\begin{definition}
  $\C$ has the \emph{amalgamation property for strong embeddings} (\emph{APS}) if for every structures $\X,\Y,\Z\in\C$ and strong embeddings $f_1:\X\to\Y$ and $f_2:\X\to\Z$, there exist $\W\in\C$ and strong embeddings $g_1:\Y\to\W$, $g_2:\Z\to\W$, such that $g_1\circ f_1=g_2\circ f_2$.
  \begin{center}\begin{tikzcd}
    \X \arrow[d, "f_1"'] \arrow[r, "f_2"] & \Z \arrow[d, "g_2", dashed]\\
    \Y \arrow[r, "g_1", dashed] & \W
  \end{tikzcd}\end{center}
\end{definition}

\begin{fact}
  \label{fact:rs}
  Suppose $\C$ satisfies the following conditions:
  \begin{enumerate}
    \item $\C$ contains an initial structure, i.e., a structure $\I\in\C^{\fin}$ such that $\I\leq\X$ for every $\X\in\C$;
    \item $\C$ has the amalgamation property for strong embeddings;
    \item $\C$ is closed under chains of strong embeddings.
  \end{enumerate}
  Then $\C$ contains a rich structure.
\end{fact}

\subsection{Facts in algebra}

In this subsection, we introduce useful tools from algebra. We are interested in the following topics: the Conjecture on Intersections with Tori (CIT), the Mordell--Lang property, the Schanuel Conjecture, and the Ax--Schanuel theorem.

We work in an algebraically closed field of characteristic $0$, denoted by $K$. Consider the algebraic $n$-torus, i.e., $(K^\times)^n$. The following proposition gives a characterization of their algebraic subgroups.

For a $(k\times n)$-integer matrix $M=(m_{ij})$, define a map from the $n$-torus to the $k$-torus:
\begin{align*}
  M(-):\;\;\;\;(K^\times)^n\;\;&\longrightarrow\;\;(K^\times)^k\\
  (x_1,...,x_n)&\longmapsto\left(\prod_{i=1}^n x_i^{m_{1i}},...,\prod_{i=1}^n x_i^{m_{ki}}\right).
\end{align*}

\begin{fact}[{\cite[Section 3.2]{Bo06}}]
  \label{fact:rank}
  Every algebraic subgroup $A$ of the $n$-torus $(K^\times)^n$ is defined by a system of equations written as follows:
  \[M(\bx)=\mathds{1},\]
  where $\bx$ is an $n$-tuple of variables and $M$ is a $(k\times n)$-integer matrix, for some $k\leq n$. Moreover, $M$ can be taken to have full rank, and then the dimension of $A$ equals $n-k$.
\end{fact}

Let $M_A$ denote any full-rank integer matrix such that $A$ is defined by equations $M_A(\bx)=\mathds{1}$.

The Conjecture on Intersections with Tori is a special case of the Zilber--Pink conjecture. It predicts that atypical components of intersections are controlled by finitely many proper algebraic subgroups. It was proposed by Zilber in 2002 \cite[Conjecture 1]{Zi02}.

\begin{conjecture}[CIT]
  Given an algebraic variety $W\se K^n$ defined over $\bQ$, there is a finite collection $\mu(W)$ of proper algebraic subgroups of $(K^\times)^n$, satisfying the following property:

  If $S$ is an atypical component of the intersection of $W$ and a proper algebraic subgroup $A\se(K^\times)^n$, i.e.,
  \[\dim^K S>\dim^K W+\dim^K A-n,\]
  then $S$ is contained in some $B\in\mu(W)$.
\end{conjecture}

Yet CIT is still open. Nevertheless, the following weaker form always holds.

\begin{fact}[Weak CIT, {\cite[Theorem 4.6]{Ki09}}]
  \label{fact:WCIT}
  Given a constructible family $(W_{\by})_{\by\in P}$ in $K^n$, i.e., constructible $W\se K^{n+m}$ and $P\se K^m$, there is a finite collection $\mu(W)$ of proper algebraic subgroups of $(K^\times)^n$, satisfying the following property:

  For every $\bc\in P$ and every proper algebraic subgroup $A\se(K^\times)^n$, if $S$ is an atypical component of the intersection of $W_{\bc}$ and a coset $\alpha A$, then there exist $B\in\mu(W)$ and a constant $\beta$, such that $S$ is contained in $\beta B$. Moreover, $S$ is typical with respect to $\beta B$.
\end{fact}

Note that Kirby's result in \cite{Ki09} is for a constructible family of varieties. To verify it for a constructible family, one checks that the uniform Schanuel property (\cite[Theorem 4.3]{Ki09}) holds for a constructible family, and the weak CIT is a direct consequence.

The Mordell--Lang property plays an important role in solving finite rank subgroups of an algebraic group. We are interested in the case of characteristic zero, which is exactly the following theorem proven by Laurent.

\begin{fact}[Mordell--Lang, {\cite{La83}}]
  \label{fact:ml}
  Let $\Gamma$ be a finite rank subgroup of $(K^\times)^n$ and let $W$ be a subvariety of $(K^\times)^n$. Then there exist a natural number $r$, elements $\gamma_1,...,\gamma_r$ of $\Gamma$, and algebraic subgroups $A_1,...,A_r$ of $(K^\times)^n$, such that $\gamma_i A_i\se W$ and
  \[W\cap\Gamma=\bigcup_{i=1}^r\gamma_i(A_i\cap\Gamma).\]
\end{fact}

The Schanuel Conjecture, if proved, would settle many well-known problems in transcendental number theory. Its proof is, however, beyond the scope of this paper. The following weaker form of the conjecture is used in the proof of Lemma \ref{le:lb}.

\begin{conjecture}[SC$_K$]
  \label{con:SCK}
  Let $K$ be a subfield of $\bC$ of finite transcendence degree. Then for all $\bx\se\bC$, we have
  \[\trd(e^{\bx})+\ld_K(\bx)-\ld_\bQ(\bx)\geq-\trd(K).\]
\end{conjecture}

By Bays, Kirby, and Wilkie \cite[Theorem 1.3]{Ki10}, the above property holds in the following special case.

\begin{fact}
  SC$_K$ holds for $K=\bQ(\beta i)$ if $\beta$ is exponentially transcendental.
\end{fact}

The Ax--Schanuel theorem gives a functional analogue of the Schanuel Conjecture. We say functions $f_1,...,f_n\in\bC[[x_1,...,x_m]]$ are \emph{$\bQ$-linearly dependent modulo $\bC$} if some of their $\bQ$-linear combination lies in $\bC$.

\begin{fact}[Ax--Schanuel, \cite{Ax71}]
  \label{fact:ax}
  Let $f_1,...,f_n\in\bC[[x_1,...,x_m]]$ be $\bQ$-linearly independent modulo $\bC$. Then
  \[\trd_\bC(f_1,...,f_n,e^{f_1},...,e^{f_n})\geq n+\rk J(f_1,...,f_n).\]
\end{fact}

\section{Construction}
\label{ch:cons}

In this section, we apply the Hrushovski construction to expansions of real closed fields by divisible multiplicative subgroups. With a suitable choice of predimension function, we establish the existence of a rich structure. Following Poizat's convention, we call it green points in the reals.

\subsection{Setup}
\label{sec:st}

Consider a structure $(\R,G)$ where $\R=(R,<,+,-,\cdot,0,1)$ is a real closed field, and $G$ is a binary predicate. Let us identify $R^2$ with an algebraically closed field, denoted by $K$. The identity $\mathds{1}$ is $(1,0)$. Suppose that $G$ defines a divisible subgroup of $K^\times$. Let $T_G$ be the theory stating that $\R\vDash\rcf$ and $G$ is a divisible group with respect to complex multiplication. Let $\C$ be the class of all such structures:
\[\C=\{\,(\R,G)\mid(\R,G)\vDash T_G\,\}.\]
Let $\cl$ be the divisible hull on $(R^2)^\times$.

Compared to the complex case, the real case contains more information. Most importantly, complex conjugation becomes definable. Therefore, the predimension function must include the algebraic information of both real and imaginary axes.

For any finite tuple $\ba=(a_1,...,a_n)=\bigl((x_1,y_1),...,(x_n,y_n)\bigr)$ of $R^2$, let $\trd(\ba)=\trd(x_1,y_1,...,x_n,y_n)$. Take the following predimension function:
\[\delta(\ba)=\trd(\ba)-\md(\cl(\ba)\cap G).\]
Clearly, $\delta$ is well-defined with respect to $\cl$. Also, $\delta$ is submodular for $\cl$-closed sets, because $\trd$ is submodular and $\md$ is modular. Thus, $\delta$ is indeed a predimension function with respect to $\cl$.

In order to apply the Hrushovski construction, we need to restrict to a class of structures where the predimension function is non-negative. Define the following class of structures:
\begin{align*}
  &\C_0=\{\,(\R,G)\in\C\mid G\text{ is torsion-free and }\delta(\ba)\geq 0\text{ for all }\ba\se(R^2)^\times\,\},\\
  &\Cf_0=\{\,(\R,G)\in\C_0\mid\R\text{ is finitely rcl-generated}\,\}.
\end{align*}

\begin{lemma}
  \label{le:id}
  Let $\X$ be a structure in $\C_0$. Then $G^\X\cap\ag{\bQ}^2=\{\mathds{1}\}$.
\end{lemma}
\begin{proof}
  Any torsion-free element of $(R^2)^\times$ has multiplicative degree $1$. Therefore, if $z$ is a non-trivial element of $G^\X\cap\ag{\bQ}^2$, then $\delta(z)=-1$. Contradicting $\X\in\C_0$.
\end{proof}

Consequently, $\C_0$ can also be written as $\{\,\X\in\C\mid\I\leq\X\,\}$ where
\[\I=(\ag{\bQ},\{\mathds{1}\})\in\Cf_0.\]

For any structures $\X\se\Y$ in $\C_0$ and any finite tuple $\ba$ in $Y^2$, define
\[\delta(\ba/\X)=\trd(\ba/X)-\md(\cl(\ba)\cap G^\Y/G^\X).\]

For any structures $\X\se\Y\in\C_0$, we say that $\Y$ is a strong extension of $\X$ if for all $\ba\se Y^2$ we have $\delta(\ba/\X)\geq 0$. This is compatible with Definition \ref{def:ss}. For $\X\se\Y$ to be a strong extension, it suffices to check that $\delta(\ba/\X)\geq 0$ for all multiplicatively independent $\ba\se G^\Y$.

In order to classify strong extensions, we also define $\delta$ for structures in $\Cf_0$. For any structures $\X\leq\Y$ in $\Cf_0$, let
\[\delta(\Y/\X)=\trd(Y/X)-\md(G^\Y/G^\X).\]
Since $\Y$ is finitely generated, we see that $\delta(\Y/\X)$ is finite.

\begin{lemma}
  \label{le:hull0}
  Let $\A\in\C_0$ and $X\se A$. Then we have $G^\A\cap\ag{[X]^\A}^2=G^\A\cap[X]^\A$ and $(\ag{[X]^\A},G^\ag{[X]^\A})\leq\A$.
\end{lemma}
\begin{proof}
  Because $[X]^\A$ is strong in $\A$, we have $\trd(\ba/[X]^\A)\geq\md(\ba/G^\A\cap[X]^\A)$ for every $\ba\se G^\A$. The assertion follows.
\end{proof}

\subsection{Existence of a rich structure}

Our goal is to prove the existence of a rich structure in $\C_0$. We want to apply Fact \ref{fact:rs}. By Lemma \ref{le:id}, $\Cf_0$ contains an initial structure. Also, $\C_0$ is closed under chains of strong embeddings (Fact \ref{fact:str}). To obtain a rich structure, it suffices to show that $\C_0$ has APS.

\begin{definition}
  A theory $T$ has the \emph{free amalgamation property} if for every structures $\X,\Y,\Z\vDash T$ and embeddings $f_1:\X\to\Y$, $f_2:\X\to\Z$, there exist $\W\vDash T$ and embeddings $g_1:\Y\to\W$, $g_2:\Z\to\W$, such that
  \begin{itemize}
    \item $g_1\circ f_1=g_2\circ f_2$, and
    \item $g_1(Y)$ and $g_2(Z)$ are free over $g_1\circ f_1(X)$, i.e., $g_1(Y)\forkindep^{\operatorname{acl}}_{g_1\circ f_1(X)}g_2(Z)$.
  \end{itemize}
\end{definition}

Since RCF is o-minimal, it has the free amalgamation property \cite[Lemma 1.2]{Pi88}. We are now ready to prove APS for $\C_0$.

\begin{proposition}
  \label{prop:am}
  $\C_0$ has the amalgamation property for strong embeddings.
\end{proposition}
\begin{proof}
  Let $(\X,G^\X),(\Y,G^\Y),(\Z,G^\Z)$ be structures in $\C_0$ with $(\X,G^\X)$ strongly embedded into $(\Y,G^\Y)$ and $(\Z,G^\Z)$. Let $\W$ be a free amalgam of $\Y$ and $\Z$ over $\X$. Identify $X,Y,Z$ with their images in $W$. Define $G^\W=G^\Y\cdot G^\Z$.

  By freeness, $(\Y,G^\Y)$ and $(\Z,G^\Z)$ are embedded into $(\W,G^\W)$. Let us show that the embedding $\Z\to\W$ is strong. Strongness of $\Y\to\W$ follows similarly.

  Let $w_1,...,w_n\in G^\W$. There exist $y_1,...,y_n\in G^\Y$ and $z_1,...,z_n\in G^\Z$ such that $w_i=y_i\cdot z_i$. Thus, $\delta(\bw/\Z)= \delta(\by/\Z)=\delta(\by/\X)$, which is non-negative because $\X\to\Y$ is strong.
\end{proof}

Applying Fact \ref{fact:rs}, we arrive at the conclusion.

\begin{corollary}
  $\C_0$ contains a rich structure.
\end{corollary}

The following lemma is a by-product of free amalgamation.

\begin{lemma}
  \label{le:inf}
  Let $\U\in\C_0$ be a rich structure, and let $\X,\Y\leq\U$ be structures in $\Cf_0$ with $\X\leq\Y$. Then there is an infinite collection of copies $(\Y_i\leq\U)_{i\in\omega}$ such that $\X\leq\Y_i$, $\Y\cong\Y_i$, and $\bigcup_{j<i}Y_j\forkindep^{\rcl}_X Y_i$ for all $i<\omega$.
\end{lemma}
\begin{proof}
  We inductively construct a chain $\W_0\leq\W_1\leq...\leq\U$ and strong substructures $(\Y_i\leq\W_i)_{i\in\omega}$.
  \begin{itemize}
    \item $\W_0=\Y_0:=\Y$.
    \item Suppose we already have $\W_{n-1}\leq\U$. By Proposition \ref{prop:am}, there is an amalgam $\W_n'$ of $\W_{n-1}$ and $\Y$ over $\X$, and strong embeddings $f_n:\W_{n-1}\to\W_n'$, $g_n:\Y\to\W_n'$ satisfying $f_n(W_{n-1})\forkindep^{\rcl}_X g_n(Y)$. Richness of $\U$ implies that there is a strong embedding $h_n:\W_n'\to\U$ extending $\W_{n-1}\leq\U$. Take $\W_n:=h_n(\W_n')$ and $\Y_n:=h_n(g_n(\Y))$.\qedhere
  \end{itemize}
\end{proof}

\subsection{Localization}\label{sec:lo}

Recall that structures in $\C_0$ must satisfy $\delta(\ba)\geq 0$ for every finite tuple $\ba$, which may be too strong for some proposed structures. As a solution, we may construct rich structures based on a fixed small structure, which possibly carries torsion information.

More precisely, recall that $(\R,G)\vDash T_G$ if and only if $\R\vDash\rcf$ and $G$ is a divisible subgroup of $(R^2)^\times$. Fix a structure $(\B,G^\B)$ of $T_G$. Expand the language by constants $c_\B=\{c_b\}_{b\in B}$, and consider $L(c_\B)$-structures extending $\B$, i.e.,
\[\C_\B'=\{\,(\R,G^\R)\vDash T_G\mid(\R,G^\R)\text{ extends }(\B,G^\B)\,\}.\]
Clearly, this is an elementary class.

We construct rich structures with respect to $\delta_\B$ given by
\[\delta_\B(\ba)=\trd(\ba/B)-\md(\cl(\ba)\cap G^\R/G^\B).\]
Define the following class of structures.
\begin{align*}
  \C_\B&=\{\,(\R,G^\R)\in\C_\B'\mid\delta_\B(\ba)\geq 0\text{, for all }\ba\se R^2\,\}\\
  \Cf_\B&=\{\,(\R,G^\R)\in\C_\B\mid\R\text{ is finitely rcl-generated over }\B\,\}
\end{align*}

The construction of a rich structure is analogous to the previous subsection, with $\C_0$ replaced by $\C_\B$.
\begin{corollary}
  $\C_\B$ contains a rich structure.
\end{corollary}

\section{Rotundity}\label{ch:rt}

In this section, we introduce rotundity, an algebraic property used to characterize the type of strong extensions. Let $\R=(R,<,+,-,\cdot,0,1)$ be a real closed field, and let $K=R^2$ be the corresponding algebraically closed field.

For a $(k\times n)$-integer matrix $M=(m_{ij})$, define the following map:
\begin{align*}
  M(-):\;\;\;\;\;\;\;\;R^{2n}\;(=K^n)\;\;\;&\longrightarrow\;\;\;(K^k=)\;R^{2k}\\
  (x_1,y_1,...,x_n,y_n)&\longmapsto\left(\prod_{j=1}^n(x_j+iy_j)^{m_{1j}},...,\prod_{j=1}^n(x_j+iy_j)^{m_{kj}}\right).
\end{align*}

\begin{definition}
  \label{de:rtd}
  A semialgebraic set $V\se R^{2n}$ is \emph{rotund} if $\dim^\R(M(V))\geq k$ for every $M\in\bZ^{k\times n}$ of rank $k$.
\end{definition}

Take a tuple $\ba\in R^{2n}$ and a parameter set $X\se R$. Let $V\se R^{2n}$ be the algebraic locus of $\ba$ over $X$. Then $V$ is rotund if and only if
\[\trd(M(\ba)/X)\geq k\text{ for every }M\in\bZ^{k\times n}\text{ of rank }k.\]

\begin{proposition}
  \label{prop:rtd}
  Let $\X\se\Y$ be structures in $\C_\B$ such that $Y=\ag{X,\ba}$, where $\ba\in (G^\Y)^n$ is a $\cl$-basis of $G^\Y$ over $G^\X$. Let $V$ be the algebraic locus of $\ba$ over $X$. Then $V$ is rotund if and only if $\X\leq\Y$.
\end{proposition}
\begin{proof}
  $(\Rightarrow)$ For every $\bb\se G^\Y$, there is some $M\in\bZ^{k\times n}$ of rank $k$ satisfying
  \[G^\X\cdot\cl(\bb)=G^\X\cdot\cl(M(\ba)).\]

  Since $\ba$ is multiplicatively independent over $G^\X$, we have
  \[\md(M(\ba)/G^\X)=k.\]

  Rotundity of $V$ implies
  \[\delta(\bb/\X)=\delta(M(\ba)/\X)=\trd(M(\ba)/X)-k\geq 0.\]

  $(\Leftarrow)$ For every $M$ of rank $k$, check that
  \[\trd(M(\ba)/X)\geq\md(M(\ba)/G^\X)=k.\qedhere\]
\end{proof}

On the other hand, take a rotund irreducible algebraic set $V\se R^{2n}$ defined over $X$. Then every generic point $\ba$ of $V$ is multiplicatively independent over $X^2$. Indeed, if $M(\ba)\in X^2$ for some $M$, then $\trd(M(\ba)/X)=0$ contradicting rotundity.

\subsection{Variety and transcendence}

Consider a point $\ba=((x_1,y_1),...,(x_n,y_n))\in R^{2n}$. Denote
\[\ha:=\bx+i\by=(x_1+iy_1,...,x_n+iy_n)\in K^n.\]
Clearly, $\trd^\R(\ba)=\trd^K(\ha\ha^c)$. We hope to infer the algebraic locus of $\ha$ and $\ha^c$ from that of $\ba$.

\begin{definition}
  Let $V(x_1,y_1,...,x_n,y_n)\se R^{2n}$ be an algebraic set. After identifying $R^2$ with an algebraically closed field $K$, we can view $V$ as a subset of $K^n$, i.e., the set
  \[\{\,(x_1+iy_1,...,x_n+iy_n)\in K^n\mid(x_1,y_1,...,x_n,y_n)\in V\,\}.\]
  Define $\widehat{V}\se K^n$ to be the Zariski closure of this set.
\end{definition}

\begin{definition}
  \label{def:tilde}
  Let $V(x_1,y_1,...,x_n,y_n)\se R^{2n}$ be an algebraic set. Define $\widetilde{V}\se K^{2n}$ to be the Zariski closure of the following set:
  \[\{\,(\bx+i\by,\bx-i\by)\in K^{2n}\mid(x_1,y_1,...,x_n,y_n)\in V\,\}.\]
\end{definition}

Let $V\se R^{2n}$ be an irreducible algebraic set. Then $\widehat{V}$ and $\widetilde{V}$ are also irreducible. Furthermore, if $(x_1',y_1',...,x_n',y_n')$ is a generic point of $V$, then $(\bx'+i\by')$ and $(\bx'+i\by',\bx'-i\by')$ are generic in $\widehat{V}$ and $\widetilde{V}$ respectively.

To acquire irreducible algebraic sets from a semialgebraic set, we use analytic cell decomposition and properties of Nash manifolds.

\begin{lemma}\label{le:dimNash}
  Let $V\se R^n$ be a semialgebraic set that is also a semialgebraically connected analytic manifold. Then the Zariski closure $V^{\zar}$ of $V$ in $R^n$ is irreducible, and $V^{\zar}$ has the same dimension as $V$.
\end{lemma}
\begin{proof}
  $V$ is a semialgebraic analytic manifold, i.e., a Nash manifold. The assertion follows from \cite[Proposition 8.4.1]{Bo98}.
\end{proof}

As a remark, such $V$ is rotund if and only if $V^{\zar}$ is rotund.

\begin{definition}
  A \emph{rotund block} is a rotund semialgebraic set $V\se R^{2n}$ that is also a semialgebraically connected analytic manifold of dimension $n$.
\end{definition}

\subsection{Definability of rotundity}

To show that rotundity is a definable property, the challenge is that there are infinitely many matrices. Our strategy is to apply the weak CIT.

Let $V_{\bt}(x_1,y_1,...,x_n,y_n)\se R^{2n}$ be an algebraic set where $\bt\se R^m$ denotes a tuple of coefficients. Then $(\widetilde{V}_{\bt})_{\bt\in K^m}$ is a constructible family in $K^{2n}$. Applying Fact \ref{fact:WCIT} to $(\widetilde{V}_{\bt\bz}=\widetilde{V}_{\bt}(-,\bz))_{\bt\bz\in K^{m+n}}$, we obtain a finite set $\mu(\widetilde{V})$ of proper algebraic subgroups, independent of the choices of $\bt\bz$. Recall Fact \ref{fact:rank}: for each $B\in\mu(\widetilde{V})$, we write $M_B$ for a matrix defining $B$. Let $\theta_{V}(\bt)$ be the conjunction of the following statements:
\begin{itemize}
  \item $V_{\bt}$ is irreducible.
  \item For each $B\in\mu(\widetilde{V})$, the dimension of $M_B(V_{\bt})$ is at least $\rk(M_B)$.
\end{itemize}
\begin{lemma}\label{le:defrtd0}
  Let $V(x_1,y_1,...,x_n,y_n;t_1,...,t_m)$ be an algebraic set (defined over the empty set). For all $\R\vDash\rcf$ and $\bc\in R^m$, we have
  \begin{center}
    $V_{\bc}$ is irreducible and rotund if and only if $\R\vDash\theta_V(\bc)$.
  \end{center}
\end{lemma}
\begin{proof}
  $(\Rightarrow)$ is clear.

  $(\Leftarrow)$ Assume $\theta_V(\bc)$ holds. Let us omit the coefficients $\bc$ and simply write $V,\widetilde{V}$ for $V_{\bc},\widetilde{V}_{\bc}$.
  
  Towards a contradiction, suppose $V$ is not rotund. Then there is a full-rank matrix $M$ such that $\dim^\R(M(V))<\rk(M)$. Let $A\se(K^\times)^n$ be the proper algebraic subgroup defined by $M(\bz)=\mathds{1}$.

  Let $\ba$ be generic in $V$ over $\bc$. Then $\ha$ is generic in $\widetilde{V}_{\ha^c}$ over $\bc\ha^c$. Let $S$ be an irreducible component of $\widetilde{V}_{\ha^c}\cap\ha A$ containing $\ha$. By the theorem of fiber dimension \cite[Theorem 1.25]{Sh13}, we have the following equality:
  \[\dim^K(S)=\dim^K(\widetilde{V}_{\ha^c}\cap\ha A)=\dim^K(\widetilde{V}_{\ha^c})-\dim^K(M\widetilde{V}_{\ha^c}).\]
  For every ($m\times n$)-matrix $\mathbf{M}$ of rank $m$, we have
  \begin{equation}\label{eq:dimMV}
    \dim^\R(\mathbf{M}V)=\trd_{\bc}^\R(\mathbf{M}\ba)=\trd_{\bc}^K(\mathbf{M}\ha,\mathbf{M}\ha^c).
  \end{equation}
  Thus,
  \[n-\dim^K(\ha A)=\rk(M)>\dim^\R(M V)\geq\trd_{\bc}^K(M\ha/\ha^c)=\dim^K(M\widetilde{V}_{\ha^c}).\]
  The following calculation shows that $S$ is an atypical component of the intersection of $\widetilde{V}_{\ha^c}$ and $\ha A$:
  \[\dim^K(S)=\dim^K(\widetilde{V}_{\ha^c})-\dim^K(M\widetilde{V}_{\ha^c})>\dim^K(\widetilde{V}_{\ha^c})+\dim^K(\ha A)-n.\]
  By Fact \ref{fact:WCIT}, there exists $B\in\mu(\widetilde{V})$ and some $\hb$ such that $S$ is contained in $\hb B$.
  
  For every coset of algebraic subgroup $\ha'A'\se\ha A$ such that $S\se\ha'A'$, we have
  \[\dim^K(\widetilde{V}_{\ha^c}\cap\ha A)\geq\dim^K(\widetilde{V}_{\ha^c}\cap\ha'A')\geq\dim^K(S)=\dim^K(\widetilde{V}_{\ha^c}\cap\ha A).\]
  Therefore, we may assume $A$ to be the smallest algebraic subgroup having a coset containing $S$. Thus, $\ha A\se\hb B$, and $\ha B=\hb B$. The typicality of $S$ with respect to $\ha B$ gives the following equality:
  \begin{equation}\label{eq:dimVcap}
    \dim^K(\widetilde{V}_{\ha^c}\cap\ha A)=\dim^K(S)=\dim^K(\widetilde{V}_{\ha^c}\cap\ha B)+\dim^K(A)-\dim^K(B).
  \end{equation}
  
  Suppose $B$ is defined by $N(\bz)=\mathds{1}$. Since $\ker(M)\se\ker(N)$, we have that $N$ factors through $M$, i.e., $N=N'\circ M$ for some $N'$. Hence,
  {\allowdisplaybreaks\begin{align*}
    \rk(M)-\rk(N)&=\dim^K(B)-\dim^K(A)\\
    &=\dim^K(\widetilde{V}_{\ha^c}\cap\ha B)-\dim^K(\widetilde{V}_{\ha^c}\cap\ha A)&&\text{by \eqref{eq:dimVcap}}\\
    &=\dim^K(M\widetilde{V}_{\ha^c})-\dim^K(N\widetilde{V}_{\ha^c})&&\text{by fiber dimension}\\
    &=\trd_{\bc}^K(M\ha/\ha^c)-\trd_{\bc}^K(N\ha/\ha^c)&&\text{by genericity}\\
    &=\trd_{\bc}^K(M\ha/N\ha,\ha^c)\\
    &\leq\trd_{\bc}^K(M\ha/N\ha,N\ha^c)\\
    &\leq\trd_{\bc}^K(M\ha/N\ha)\\
    &\leq\md(M\ha/N\ha)\\
    &\leq\rk(M)-\rk(N),
  \end{align*}}
  and hence all terms in the inequality are equal. In particular,
  \begin{equation}\label{eq:rkM}
    \rk(M)-\rk(N)=\trd_{\bc}^K(M\ha/N\ha,N\ha^c)=\trd_{\bc}^K(M\ha,N\ha^c)-\trd_{\bc}^K(N\ha,N\ha^c).
  \end{equation}
  Put everything together, we get
  \begin{align*}
    \dim^\R(M(V))-\dim^\R(N(V))&=\trd_{\bc}^K(M\ha,M\ha^c)-\trd_{\bc}^K(N\ha,N\ha^c)&&\text{by \eqref{eq:dimMV}}\\
    &\geq\trd_{\bc}^K(M\ha,N\ha^c)-\trd_{\bc}^K(N\ha,N\ha^c)\\
    &=\rk(M)-\rk(N).&&\text{by \eqref{eq:rkM}}
  \end{align*}
  As a result,
  \[\dim^\R(N(V))-\rk(N)\leq\dim^\R(M(V))-\rk(M)<0.\]
  However, $\theta_V(\bc)$ states that $N(V)$ has dimension at least $\rk(N)$. Contradiction.
\end{proof}

\begin{lemma}\label{le:defrtd}
  Let $W(x_1,y_1,...,x_n,y_n;w_1,...,w_k)$ be a semialgebraic set (defined over the empty set). Let $\Theta^W_V(\bw,\bt)$ be the conjunction of the following formulas:
  \begin{itemize}
    \item $V_{\bt}$ is irreducible and rotund, i.e., $\theta_V(\bt)$.
    \item $W_{\bw}$ is a semialgebraically connected analytic manifold.
    \item $W_{\bw}\se V_{\bt}\;\bigwedge\;\dim^\R W_{\bw}=\dim^\R V_{\bt}=n$.
  \end{itemize}
  Then for all $\R\vDash\rcf$ and $\bb\in R^k$, we have
  \begin{center}
    $W_{\bb}$ is a rotund block if and only if $\R\vDash\Theta^W_V(\bb,\bc)$ for some $V_{\bc}\se R^{2n}$.
  \end{center}
\end{lemma}
\begin{proof}
  Take $V_{\bc}$ to be the Zariski closure of $W_{\bb}$. The statement follows from Lemma \ref{le:dimNash}.
\end{proof}

\subsection{Generic hyperplanes}

This subsection is devoted to the existence of rotund blocks. Typical examples are generic hyperplanes defined in \cite[Section 4.2.2]{Ca11}. Let $H_{n,m}(\bx;\by_1,...,\by_m)$ denote the subvariety of $K^{n+mn}$ defined by equations $M(\bx)=\mathds{1}$, where $M$ is the $(m\times n)$-matrix with rows $\by_1,...,\by_m$.

\begin{fact}[{\cite[Remark 4.2.23]{Ca11}}]
  \label{fact:gh}
  For every algebraically closed subfield $L$ of $K$, there exist $\bc_1,...,\bc_n\in L^{2n}$ such that the following holds:
  \begin{itemize}
    \item $\dim^K(H_{2n,n}(-;\bc_1,...,\bc_n))=n$,
    \item $\dim^K(M(H_{2n,n}(-;\bc_1,...,\bc_n)))>\frac{k}{2}$, for every $M\in\bZ^{k\times 2n}$ of rank $k\in[1,2n)$.
  \end{itemize}
\end{fact}

Let $H_{n,m}'$ denote the image of $H_{n,m}\se K^{n+mn}$ in $(R^2)^{n+mn}$.

\begin{proposition}\label{prop:gen}
  For every parameter set $X$, there is a rotund block $V\se R^{4n}$ defined over $X$ such that
  \begin{itemize}
    \item $\dim^\R(V)=2n$,
    \item $\dim^\R(M(V))>k$, for every $M\in\bZ^{k\times 2n}$ of rank $k\in[1,2n)$,
    \item the projection onto its first coordinate has image $R$.
  \end{itemize}
\end{proposition}
\begin{proof}
  Let $H_{2n,n}(-;\bc_1,...,\bc_n)\se K^{2n}$ be an instance of Fact \ref{fact:gh} with $L=\ag{X}^2$. Thus,
  \begin{itemize}
    \item $\dim^\R(H_{2n,n}'(-;\bc_1,...,\bc_n))=2\dim^K(H_{2n,n}(-;\bc_1,...,\bc_n))=2n$, and
    \item for every $M\in\bZ^{k\times 2n}$ of rank $k\in[1,2n)$, we have
    \[\dim^\R(M(H_{2n,n}'(-;\bc_1,...,\bc_n)))=2\dim^K(M(H_{2n,n}(-;\bc_1,...,\bc_n)))>k.\]
  \end{itemize}
  Therefore, $H_{2n,n}'(-;\bc_1,...,\bc_n)\se R^{4n}$ is a rotund block. It is clear that the projection onto one of its coordinates has image $R$. After swapping the coordinates we obtain the desired result.
\end{proof}

\section{Rich structures}
\label{ch:rs}

In this section, we analyze basic properties of strong extensions and rich structures. Although the following statements are presented for rich structures in $\C_0$, the same also holds for $\C_\B$ with minor modification of the proofs.

\subsection{Minimal extension}

Every proper strong extension in $\Cf_0$ can be decomposed into a finite sequence of minimal strong extensions.

\begin{proposition}
  \label{prop:min}
  Suppose $\X\leq\Y\in\Cf_0$ is minimal. Then exactly one of the following holds:
  \begin{enumerate}
    \item $\delta(\Y/\X)=0$, and there exists a $\cl$-basis $\ba\in (G^\Y)^n$ of $G^\Y$ over $G^\X$, such that $\trd(\ba/X)=n$ and $Y=\ag{X,\ba}$.
    \item $\delta(\Y/\X)=1$, $G^\Y=G^\X$, and there exists $y\in Y\setminus X$ such that $Y=\ag{X,y}$.
  \end{enumerate}
\end{proposition}
\begin{proof}
  $\delta(\Y/\X)$ can be either zero or positive.
  \begin{enumerate}
    \item $\delta(\Y/\X)=0$.
    
    Take a $\cl$-basis $\ba\in(G^\Y)^n$ of $G^\Y$ over $G^\X$. Observe that
    \[n=\md(\ba/G^\X)\leq\trd(\ba/X)\leq\trd(Y/X)=\md(G^\Y/G^\X)=n.\]
    Therefore, $\md(\ba/G^\X)=\trd(\ba/X)=n$.

    \item $\delta(\Y/\X)>0$.
    
    For every proper extension $\X\leq\Z\se\Y$, we must have $\delta(\Z/\X)>0$. Otherwise, if $\delta(\Z/\X)=0$, then $\X\leq\Z\leq\Y$ gives a chain of proper extensions, contradicting minimality. Thus, $\delta(\ba/\X)>0$ for all $\ba\se Y^2\setminus X^2$.

    Let $y_1,...,y_n$ be a transcendence basis of $Y$ over $X$. Then we have
    \begin{align*}
      &\delta(\ag{X,y_1}/\X)=1,\\
      &\delta(\ba/\ag{X,y_1})\geq 0,\text{ for every }\ba\se Y^2.
    \end{align*}
    Hence, $\left(\ag{X,y_1},G^\ag{X,y_1}\right)\leq\Y$. We must have $Y=\ag{X,y_1}$, or it would violate minimality. To see $G^\Y=G^\X$, notice that
    \[1=\trd(Y/X)>\md(G^\Y/G^\X).\qedhere\]
  \end{enumerate}
\end{proof}

We call a minimal extension $\X\leq\Y\in\Cf_0$ \emph{prealgebraic} if $\delta(\Y/\X)=0$, and \emph{purely transcendental} if $\delta(\Y/\X)=1$.

Fix a rich structure $\U\in\C_0$.

\begin{lemma}\label{le:inters}
  Let $\Y,\Z\leq\U$ be structures in $\Cf_0$. Denote $X=Y\cap Z$. If $\X\leq\Y$ is minimal prealgebraic, then $Y\forkindep^{\rcl}_X Z$.
\end{lemma}
\begin{proof}
  Fact \ref{fact:str} implies $\X$ is also strong in $\U$. Let $\ba\in(G^\Y)^n$ be a cl-basis of $G^\Y$ over $G^\X$ such that $\trd(\ba/X)=n$ and $Y=\ag{X,\ba}$. Since $\md$ is modular, we have
  \[\trd(\ba/Z)\geq\md(\ba/G^\Z)=\md(\ba/G^\X)=n.\qedhere\]
\end{proof}

The dimension function $d:\mathcal{P}(U^2)\to\bZ$ associated to $\delta$ is given by
\begin{align*}
d(\bx)&=\min\{\delta(\bx\bz)\mid\bz\se\U^2\}\\
&=\min\{\trd(\bx\bz)-\md(\cl(\bx\bz)\cap G^\U)\mid\bz\se\U^2\}\\
&=\min\{\trd(\bx\by)-\md(\by)\mid\by\se G^\U\}.
\end{align*}
The last equality is obtained by taking $\by$ to be a $\cl$-basis of $\cl(\bx\bz)\cap G^\U$. We can extend the domain of $d$ to $\mathcal{P}(U)$ by defining
\[d(\bx):=\min\{\trd(\bx\by)-\md(\by)\mid\by\se G^\U\},\text{ for each }\bx\se U.\]

The induced pregeometry $\Cl$ is given by the union of all prealgebraic extensions, i.e.,
\[\Cl(\bx)=\bigcup\,\{\Y\leq\U\mid\ag{[\bx]}\leq\Y\wedge\delta(\Y/\ag{[\bx]})=0\}.\]

Let $\bb\se G^\U$ and $t\in\U$. Then $t\in\U\setminus\Cl(\bb)$ if and only if
\[\trd(\ba,t/\bb)>\md(\ba/\bb),\text{ for all }\ba\se G^\U.\]

\subsection{Existential closedness}
\label{sec:ec}

In this subsection, we define a set of axioms EC in order to characterize existential closedness of rich structures.

\begin{definition}
  A structure $\U$ in $\C_0$ is said to be \emph{existentially closed with respect to strong extensions} if for every quantifier-free formula $\phi(\bx)$ and every strong extension $\U'$ of $\U$, we have that $\U'\vDash\exists\bx\,\phi(\bx)$ implies $\U\vDash\exists\bx\,\phi(\bx)$.
\end{definition}

\begin{lemma}
  \label{le:ex}
  If $\U$ is rich in $\C_0$, then $\U$ is existentially closed with respect to strong extensions.
\end{lemma}
\begin{proof}
  Let $\phi(\bx)$ be a quantifier-free formula defined over $\bc\se U$, and let $\U'$ be a strong extension of $\U$. By Lemma \ref{le:hull0}, there is a finitely rcl-generated substructure $\X$ with $X=\ag{[\bc]^{\U}}$ strong in $\U$. Take $\ba\se U'$ such that $\U'\vDash\phi(\ba)$. Fact \ref{fact:str} implies $\X\leq(\ag{X,\ba},G|_\ag{X,\ba})$, and hence $(\ag{X,\ba},G|_\ag{X,\ba})$ is strongly embedded into $\U$. Take the image of $\ba$.
\end{proof}

\begin{definition}
  A structure $(\R,G)\vDash T_G$ has the \emph{EC-property} if for every rotund block $V\se R^{2n}$, the intersection $V\cap G^n$ is non-empty.
\end{definition}

\begin{lemma}
  \label{le:EC}
  If $(\R,G)$ is rich in $\C_0$, then $(\R,G)$ has the EC-property.
\end{lemma}
\begin{proof}
  Let $V$ be a rotund block over $R$. Let $\tilde{\R}$ be an extension of $\R$ containing a generic point $\ba=(a_1,...,a_n)\in(\tilde{R}^2)^n$ of $V$. We pick a sequence of compatible roots of $a_k$ for each $k=1,...,n$. More precisely, let $(a^{[i]}_k)_{i\geq 1}$ be a sequence of points in $\tilde{R}^2$ such that $a^{[1]}_k=a_k$ and $(a^{[i\cdot j]}_k)^j=a^{[i]}_k$ for all $i,j\geq 1$. Let $\R'$ be the real closed field generated by $R\cup\ba$ and let $G^{\R'}$ be the group generated by $G^\R\cup(\ba^{[i]})_{i\geq 1}$. Proposition \ref{prop:rtd} implies $(\R,G^\R)\leq(\R',G^{\R'})$. By Lemma \ref{le:ex}, we conclude that $\U\vDash\exists\bx,\,\bx\in V\cap G^n$.
\end{proof}

Lemma \ref{le:defrtd} enables us to express the EC-property in a set of first order sentences denoted by EC.

\subsection{Back and forth}
\label{sec:bf}

Rich structures are back-and-forth equivalent. For $\U_1,\U_2\in\C_0$, let $\F(\U_1,\U_2)$ be the following set of partial isomorphisms:
\[\F(\U_1,\U_2)=\{\,f:\X_1\xrightarrow{\sim}\X_2\mid\X_i\leq\U_i,\,\X_i\in\Cf_0\,\}.\]
\begin{lemma}
  \label{le:bf}
  Let $\U_1,\U_2$ be rich structures in $\C_0$. Then $\F(\U_1,\U_2)$ forms a back-and-forth system.
\end{lemma}
\begin{proof}
  Let $f:\X_1\xrightarrow{\sim}\X_2$ be a partial isomorphism in $\F(\U_1,\U_2)$. Take $a_1\in(U_1)^2$. Let us show that there exists a partial isomorphism $g\in\F(\U_1,\U_2)$ extending $f$ with $a_1\se\operatorname{Dom}(g)$. By Fact \ref{fact:min}, there is a tuple $\ba=(a_1,...,a_n)\se(U_1)^2$ such that $(\ag{X_1,\ba},G|_\ag{X_1,\ba})\leq\U_1$.

  Notice that $f^{-1}$ gives a strong embedding $X_2\to\ag{X_1,\ba}$. By richness of $\U_2$, there is a strong embedding $g:\ag{X_1,\ba}\to U_2$ such that $g\circ f^{-1}$ is the inclusion $X_2\se U_2$. Thus, $g:\ag{X_1,\ba}\xrightarrow{\sim}\ag{X_2,g(\ba)}$ provides a desired partial isomorphism.
\end{proof}

Observe that $\F(\U_1,\U_2)$ is nonempty. Indeed, it must contain the identity map on the initial structure.

\begin{corollary}
  All rich structures in $\C_0$ are elementarily equivalent.
\end{corollary}

\section{Axiomatization}
\label{ch:th}

\subsection{Axiomatizing \texorpdfstring{$\C_0$}{C0}}
\label{sec:axc0}

In this subsection, we find the theory of class $\C_0$, using the information provided by Fact \ref{fact:WCIT} (weak CIT).

Recall that for each $\ba\in G^n$, we want $\trd(\ba)\geq\md(\ba)$. The idea is that, if $\ba$ satisfies too many algebraic relations, then it must be multiplicatively dependent.

For each irreducible variety $V\se R^{2n}$ over $\bQ$, by Fact \ref{fact:WCIT} applied to $\widehat{V}$, we obtain a finite collection $\mu(\widehat{V})$ of proper algebraic subgroups. Denote
\[r_V:=\dim^\R(V)-\dim^K(\widehat{V})=\dim^K(\widetilde{V})-\dim^K(\widehat{V}).\]
Let $\phi_V$ be the following sentence:
\begin{align*}
&\forall\,(x_1,y_1,...,x_n,y_n)\in V\cap G^n,\\
&\dim^K\left(\widetilde{V}(\bx+i\by,-)\right)=r_V\,\longrightarrow\bigvee_{B\in\mu(\widehat{V})}(x_1+iy_1,...,x_n+iy_n)\in B.
\end{align*}

Let $\Phi$ consist of all $\phi_V$ for irreducible $V\se R^{2n}$ over $\bQ$ with dimension $\dim^\R(V)<n$, and let $T_0$ be the union of $T_G$ and $\Phi$.

\begin{lemma}
  \label{le:c0}
  For every $(\R,G)\in\C$,
  \[(\R,G)\vDash\Phi\text{ if and only if }\delta(\ba)\geq 0\text{ for all }\ba\se(R^2)^\times.\]
\end{lemma}
\begin{proof}
  $(\Rightarrow)$ It suffices to show that $\delta(\ba)\geq 0$ for all multiplicatively independent $\ba\in G^n$. Suppose on the contrary that there exists multiplicatively independent $\ba\in G^n$ with $\trd^\R(\ba)<n$. Let $V\se R^{2n}$ be the algebraic locus of $\ba$ over $\bQ$. Then $\dim^\R(V)=\trd^\R(\ba)<n$. Since $(\ha,\ha^c)$ and $\ha$ are generic points of $\widetilde{V}$ and $\widehat{V}$ respectively, the following equality holds:
  \[r_V=\trd^K(\ha\ha^c)-\trd^K(\ha)=\trd^K(\ha^c/\ha)=\dim^K(\widetilde{V}(\ha,-)).\]
  Then $\phi_V$ states that $\ba$ must be multiplicatively dependent. Contradiction.

  $(\Leftarrow)$ Let $V\se R^{2n}$ be an irreducible variety over $\bQ$ with $\dim^\R(V)<n$. Suppose $\ba=(x_1,y_1,...,x_n,y_n)$ belongs to $V\cap G^n$ and $\widetilde{V}(\ha,-)$ has dimension $r_V$. Since $\ha^c$ is a point in $\widetilde{V}(\ha,-)$, we have $\trd^K(\ha^c/\ha)\leq r_V$. Thus,
  \[\trd^\R(\ba)=\trd^K(\ha\ha^c)\leq\trd^K(\ha)+\dim^\R(V)-\dim^K(\widehat{V}).\]
  Because of our assumption $\delta(\ba)\geq 0$, we know
  \[\md(\ha)=\md(\ba)\leq\trd^\R(\ba)\leq\dim^\R(V)<n.\]
  Thus, $\ha$ is multiplicatively dependent. Let $A$ be an algebraic subgroup containing $\ha$ with $\dim^K(A)=\md(\ha)$. Let $S$ be an irreducible component of $\widehat{V}\cap A$ containing $\ha$. The following calculation indicates that $S$ is atypical:
  \[\dim^K S\geq\trd^K\ha\geq\dim^K\widehat{V}+\trd^\R\ba-\dim^\R V>\dim^K\widehat{V}+\dim^K A-n.\]
  By Fact \ref{fact:WCIT}, there is $B\in\mu(\widehat{V})$ and a constant $\gamma$ such that for every $\bw\in S$, we have $M_B(\bw)=\gamma$. Since $S$ is defined over $\bQ$, we see that $\gamma$ must be algebraic over $\bQ$. Moreover, $\ba\se G$ implies $\gamma\se G$. By Lemma \ref{le:id}, $\gamma=\mathds{1}$.
\end{proof}

\begin{corollary}
  $T_0$ axiomatizes the class $\C_0$.
\end{corollary}

\subsection{Axiomatizing \texorpdfstring{$\C_\B$}{CB}}\label{sec:alo}

Recall the class $\C_\B$ defined in Subsection \ref{sec:lo}. We provide an axiomatization of $\C_\B$ under finiteness assumptions on $\B$.

\begin{definition}
  A structure $\B\in\C$ is said to be \emph{finitely $G$-based} if $G^\B$ has a finite $\cl$-basis $\bb$ and $B=\rcl(\bb)$.
\end{definition}

Throughout this subsection, we assume that $\B\in\C$ is finitely $G$-based and let $\bb$ be a $\cl$-basis of $G^\B$. As a consequence, $\delta_\B(\ba)=\delta(\ba/\bb)$.

Let $V\se R^{2n}$ be an irreducible variety over $\bQ(\bb)$. By Fact \ref{fact:WCIT} applied to $\widehat{V}$, we obtain a finite collection $\mu(\widehat{V})=\{B_1,...,B_s\}$ of proper algebraic subgroups. Let $M_i$ be a full-rank integer matrix such that $B_i$ is defined by equations $M_i(\bz)=\mathds{1}$. Ideally, we want the formula $M_i(\bz)\se\cl(\bb)$, which is not expressible. Instead, with $\Gamma=\cl(\bb)^n$, we apply Fact \ref{fact:ml} to $M_i(\widehat{V})$. It yields $r_i$, $\{\gamma_{i,1},...,\gamma_{i,r_i}\}$, and $\{A_{i,1},...,A_{i,r_i}\}$.

Let $\phi^{\bb}_V$ be the following sentence, where $\bz=\bx+i\by$:
\begin{align*}
  &\forall\,(x_1,y_1,...,x_n,y_n)\in V\cap G^n,\\
  &\bigwedge_{i=1}^s\dim^K(\widehat{V}\cap\bz B_i)=\dim^K(\widehat{V})-\dim^K(M_i\widehat{V})\,\bigwedge\,\dim^K(\widetilde{V}(\bz,-))=r_V\\
  &\hspace{11em}\longrightarrow\bigvee_{i=1}^s\bigvee_{j=1}^{r_i}\left(A_{i,j}\neq(K^\times)^n\bigwedge M_i(\bz)\in\gamma_{i,j}A_{i,j}\right).
\end{align*}
Let $\Phi^{\bb}$ consist of all $\phi^{\bb}_V$ for irreducible $V\se R^{2n}$ over $\bQ(\bb)$ with dimension $\dim^\R(V)<n$. Denote the axioms of $\C_\B'$ by $T_\B'$, and take $T_\B=T_\B'\cup\Phi^{\bb}$.

\begin{lemma}
  For every $(\R,G)\in\C_\B'$,
  \[\X\vDash\Phi^{\bb}\text{ if and only if }\delta(\ba/\bb)\geq 0\text{ for all }\ba\se(R^2)^\times.\]
\end{lemma}
\begin{proof}
  $(\Rightarrow)$ The same proof as Lemma \ref{le:c0}. Notice that, if $\ha$ is a generic point of $\widehat{V}$, then the dimension of fibers satisfies the following equality:
  \[\dim^K(\widehat{V}\cap\ha B_i)=\dim^K(\widehat{V})-\dim^K(M_i\widehat{V}),\text{ for all }B_i\in\mu(\widehat{V}).\]

  $(\Leftarrow)$ Let $V\se R^{2n}$ be a variety over $\bQ(\bb)$ with $\dim^\R(V)<n$. Suppose $\ba$ is in $V\cap G^n$ satisfying the second line of $\phi^{\bb}_V$. Following the same argument as the proof of $(\Leftarrow)$ in Lemma \ref{le:c0}, we take a coset $\alpha A$ of a proper algebraic subgroup $A$ such that $\ha\in\alpha A$, $\alpha\se\cl(\bb)$, and $\dim^K(A)=\md(\ba/\bb)$. The irreducible component $S$ of $\widehat{V}\cap\alpha A$ containing $\ha$ is atypical.

  By Fact \ref{fact:WCIT}, there is $B_i\in\mu(\widehat{V})$ and a constant $\gamma$ such that for every $\bw\in S$, we have $M_i(\bw)=\gamma$. Since $S$ is defined over $\bQ(\bb)$, we see that $\gamma$ must be algebraic over $\bQ(\bb)$. Moreover, $\ba\se G$ implies $\gamma\se G$, and therefore $\gamma\se\cl(\bb)$. If $M_i(\widehat{V})$ has dimension strictly less than $\rk(M_i)$, then $A_{i,j}$ must be proper for all $j$, and Fact \ref{fact:ml} gives the desired result. Let us prove that $\dim^K(M_i\widehat{V})<\rk(M_i)$.

  $S$ is atypical in the intersection of $\widehat{V}$ and $\alpha A$, and typical with respect to $\ha B_i$, i.e.,
  \begin{align*}
    \dim^K(S)&>\dim^K(\widehat{V})+\dim^K(\alpha A)-n,\\
    \dim^K(S)&=\dim^K(\widehat{V}\cap\ha B_i)+\dim^K(\alpha A\cap\ha B_i)-\dim^K(\ha B_i).
  \end{align*}
  Combining this inequality with the second line of $\phi^{\bb}_V$, we get
  \begin{align*}
    \dim^K(M_i\widehat{V})&=\dim^K(\widehat{V})-\dim^K(\widehat{V}\cap\ha B_i)\\
    &\leq\dim^K(\widehat{V})-\dim^K(\widehat{V}\cap\ha B_i)+\dim^K(\alpha A)-\dim^K(\alpha A\cap\ha B_i)\\
    &<n-\dim^K(\ha B_i)\\
    &=\rk(M_i).\qedhere
  \end{align*}
\end{proof}

\begin{corollary}
  $T_\B$ axiomatizes the class $\C_\B$.
\end{corollary}

Following the same proof of the above lemma, we obtain the following proposition, which provides first-order characterization for $\Cl$.

\begin{proposition}
  \label{prop:tre}
  There is a partial type $\Psi(x)$ over $\bb$ such that for every $t\in R$, we have that $t$ satisfies $\Psi$ if and only if
  \begin{equation}\label{eq:tr}
    \forall\,\ba\in G^n,\,\trd(\ba,t/\bb)-\md(\ba/\bb)>0.
  \end{equation}
\end{proposition}
\begin{proof}
  For each variety $V\se R^{2n}$ over $\bQ(x,\bb)$, apply Fact \ref{fact:WCIT} and Fact \ref{fact:ml} as before, keeping $\Gamma=\cl(\bb)^n$. Let $\Psi(x)$ be the conjunction of the followings:
  \begin{itemize}
    \item $x$ is transcendental over $\bQ(\bb)$, and
    \item $\phi_V^{(x,\bb)}$ for all $V$ with $\dim^\R(V)<n$.
  \end{itemize}
  Then $\Psi(t)$ holds if and only if $\trd(t/\bb)=1$, and
  \[\forall\,\ba\in G^n,\,\trd(\ba/t,\bb)-\md(\ba/\bb)\geq 0,\]
  which is equivalent to \eqref{eq:tr}.
\end{proof}

\subsection{Axiomatizing richness}\label{sec:rich}

This subsection is devoted to the proof of Theorem \ref{thm:th}.

Define the theory $T_0^{\rich}$ as the union of $T_0$ and EC. In Subsection \ref{sec:ec}, we already see that rich structures possess the EC-property. In the following, we demonstrate how the EC-property and saturation imply richness.

\begin{lemma}
  \label{le:rich}
  Let $\U$ be an $\omega$-saturated model of $T_0^{\rich}$. Then $\U$ is a rich structure in $\C_0$.
\end{lemma}
\begin{proof}
  Let $\X\leq\Y$ in $\Cf_0$ be a minimal strong extension. Assume that $\X$ is strongly embedded into $\U$. For simplicity, let us further assume that $X\se U$. We need to show that $\Y$ is strongly embedded into $\U$. By Proposition \ref{prop:min}, there are two cases.
  \begin{enumerate}[leftmargin=*]
    \item $\delta(\Y/\X)=0$.\label{it:delta0}

    Let $\bb\in(G^\Y)^n$ be a $\cl$-basis of $G^\Y$ over $G^\X$ such that $\trd(\bb/X)=n$ and $Y=\ag{X,\bb}$. Let $(\bb^{[i]})_{i\geq 1}$ be a sequence in $G^\Y$ consisting of compatible roots of $\bb$. In other words, $\bb^{[1]}=\bb$ and $(\bb^{[i\cdot j]})^j=\bb^{[i]}$ for all $i,j\geq 1$.
    \begin{claim}
      $\qftp(\bb^{[N]}|X)$ is realized in $\U$ for every $N\geq 1$.
    \end{claim}
    \begin{pfcl}
      By $\omega$-saturation, it suffices to show that $\qftp(\bb^{[N]}|X)$ is finitely satisfiable. Let $\mu$ be a finite subset of $\qftp(\bb^{[N]}|X)$. Since $\trd(\bb^{[N]}/X)=n$, by analytic cell decomposition, there exists an analytic $n$-cell $V\se R^{2n}$ defined over $X$, such that $\mu(x_1,y_1,...,x_n,y_n)$ can be deduced from $(x_1,y_1,...,x_n,y_n)\in V\cap G^n$. By definition, analytic cells are connected analytic manifolds. Since $\X\leq\Y$, by Proposition \ref{prop:rtd}, $V$ is rotund. The EC-property then implies that $V\cap(G^\U)^n$ contains an element, which is a realization of $\mu$ in $\U$.\hfill$\blacksquare$
    \end{pfcl}
    \begin{claim}
      $\qftp((\bb^{[i]})_{i\geq 1}|X)$ is realized in $\U$.
    \end{claim}
    \begin{pfcl}
      It suffices to show that $\qftp((\bb^{[i]})_{i\leq m}|X)$ is realized in $\U$ for every $m\geq 1$. Let $N=\prod_{i=1}^m i$. Let $\ba$ be a realization of $\qftp(\bb^{[N]}|X)$ in $\U$. For every $i\leq m$, define $\ba^{[i]}=\ba^{N/i}$. Then $(\ba^{[i]})_{i\leq m}$ satisfies $\qftp((\bb^{[i]})_{i\leq m}|X)$.\hfill$\blacksquare$
    \end{pfcl}

    Let $(\ba^{[i]})_{i\geq 1}\se U^2$ be a realization of $\qftp((\bb^{[i]})_{i\geq 1}|X)$. The map $(\bb^{[i]})_{i\geq 1}\mapsto(\ba^{[i]})_{i\geq 1}$ gives an embedding from $\Y$ to $\U$ fixing $\X$. Identify $\Y$ with its image. Let us show that $\Y\leq\U$.

    For every $\bw\se G^\U$, we have
    \[\delta(\bw/\Y)=\delta(\bw\bb/\X)-\delta(\Y/\X) =\delta(\bw\bb/\X)\geq 0.\]
    The last inequality follows from $\X\leq\U$.

    \item $\delta(\Y/\X)=1$.
    
    Let $y\in Y\setminus X$. Let $\Psi(t)$ be the partial type over $X$ stating that:
    \begin{itemize}
      \item $t$ realizes $\qftp(y|X)$, and
      \item $\forall\,\ba\in G^n,\,\trd(\ba,t/X)-\md(\ba/G^\X)>0$.
    \end{itemize}
    Proposition \ref{prop:tre} shows that the second line is type-definable.
    \begin{claim}
      $\Psi$ is realized in $\U$.
    \end{claim}
    \begin{pfcl}
      It suffices to show that $\Psi$ is finitely satisfiable. Fix $N\in\bN$. Let $\psi(t)$ be the following partial type.
      \begin{enumerate}
        \item $c_1<t<c_2$, for some $c_1,c_2\in X$, and
        \item for all $n<2N$, $\forall\,\ba\in G^n,\,\trd(\ba,t/X)-\md(\ba/G^\X)>0$.
      \end{enumerate}

      By Proposition \ref{prop:gen}, there is a rotund block
      \[V(b_1,...,b_{2N})=V(x_1,y_1,...,x_{2N},y_{2N})\se R^{4N}\]
      with a generic point $\bb'=(p_1,q_1,...,p_{2N},q_{2N})$ over $X$ such that
      \begin{itemize}
        \item $p_1$ is transcendental over $X$,
        \item $c_1<p_1<c_2$, and
        \item $\trd(M(\bb')/X)>k$, for every $M\in\bZ^{k\times 2N}$ of rank $k\in[1,2N)$.
      \end{itemize}
      Define $\X'$ by
      \begin{align*}
        X'&=\ag{X,\bb'},\\
        G^{\X'}&=G^\X\cdot\langle\bb'^{[i]}\mid i\geq 1\rangle_{gp},
      \end{align*}
      where $(\bb'^{[i]})$ is a sequence of compatible roots of $\bb'$. By Proposition \ref{prop:rtd}, we have $\X\leq\X'\in\Cf_0$. Notice that $\trd(\bb'/X)=\md(\bb'/G^\X)=2N$, i.e., $\delta(\X'/\X)=0$. Case \eqref{it:delta0} implies that $\X'$ is strongly embedded into $\U$. Identify $\X'$ with its image in $\U$.

      Let us show that $\U\vDash\psi(p_1)$. Clearly, (a) holds. For (b), take any $\ba\in(G^\U)^n$, where $n<2N$. Observe that $\trd(\ba,p_1/X)=\trd(\ba/\ag{X,p_1})+1$. Also, $G^\ag{X,p_1}=G^\X$. It suffices to verify the following:
      \[\delta(\ba/\ag{X,p_1})=\trd(\ba/\ag{X,p_1})-\md(\ba/G^\X)\geq 0.\]
      Take a $\cl$-basis $\ba'$ of $\cl(\ba)\cap G^{\X'}$ over $G^\X$. It follows from submodularity that
      \[\delta_\ag{X,p_1}(\ba/\ba')\geq\delta_\ag{X,p_1}(\ba/\bb')=\delta(\ba/\X')\geq 0.\]

      \begin{claim}
        $\delta_\ag{X,p_1}(\ba')=\trd(\ba'/\ag{X,p_1})-\md(\ba'/G^\X)\geq 0$.
      \end{claim}
      \begin{pfcl}
        Suppose $\ba'$ is non-empty. Then there exists a matrix $M\in\bZ^{k\times 2N}$ of rank $k$, where $1\leq k<2N$, such that $\cl(G^\X,\ba')=\cl(G^\X,M(\bb'))$. Since $\dim(M(V))>k$, we obtain
        \[\trd(\ba',p_1/X)\geq\trd(\ba'/X)=\trd(M(\bb')/X)>k=\md(\ba'/G^\X),\]
        and hence, $\trd(\ba'/\ag{X,p_1})=\trd(\ba',p_1/X)-1\geq\md(\ba'/G^\X)$.\hfill$\blacksquare$
      \end{pfcl}

      Combining the two inequalities, we get
      \[\delta(\ba/\ag{X,p_1})=\delta_\ag{X,p_1}(\ba/\ba')+\delta_\ag{X,p_1}(\ba')\geq 0.\tag*{$\blacksquare$}\]
    \end{pfcl}

    Let $y'\in U$ be a realization of $\Psi$. The map $y\mapsto y'$ induces an embedding from $\Y$ to $\U$ fixing $\X$. It is clear that $(\ag{X,y'},G^\X)\leq\U$.\qedhere
  \end{enumerate}
\end{proof}

The following results are immediate consequences of a back-and-forth system.

\begin{corollary}
  Every rich structure in $\C_0$ is $\omega$-saturated.
\end{corollary}
\begin{proof}
  Let $\U_1$ be a rich structure in $\C_0$. Take an $\omega$-saturated model $\U_2$ of $T_0^{\rich}$. Lemma \ref{le:rich} implies $\U_2$ is rich. Applying Lemma \ref{le:bf}, we see that $\A_1$ is also $\omega$-saturated.
\end{proof}

\begin{corollary}
  $T_0^{\rich}$ is complete.
\end{corollary}
\begin{proof}
  Any model of $T_0^{\rich}$ can be elementarily extended to an $\omega$-saturated model, which is rich by Lemma \ref{le:rich}. By Lemma \ref{le:bf}, rich structures are elementarily equivalent.
\end{proof}

\begin{corollary}
  Rich structures in $\C_0$ are exactly $\omega$-saturated models of $T_0^{\rich}$.
\end{corollary}

The arguments presented in this subsection trivially extend to rich structures in the localized class $\C_\B$ for finitely $G$-based $\B\in\C$. Let $T_\B^{\rich}$ denote the union of $T_\B$ and EC. It gives a complete theory of rich structures in $\C_\B$.

\section{Tameness}\label{ch:tm}

In this section, we investigate model-theoretic properties of the constructed theory. In particular, near model completeness, o-minimal open core, and strong dependence are verified. For better readability, here we demonstrate the proof for $T_0^{\rich}$, and the same proof works for the localized theory $T_\B^{\rich}$ for finitely $G$-based $\B\in\C$.

\subsection{Near model completeness}

Theories obtained by the Hrushovski construction often possess near model completeness, a weaker form of quantifier elimination.

\begin{definition}
  A theory $T$ is called \emph{near model complete} if every formula is equivalent, modulo $T$, to a Boolean combination of existential formulas.
\end{definition}

To show near model completeness of a theory $T$ constructed by predimension method, we use the following fact \cite[Proposition 2.3.11]{Ca11}.

\begin{fact}
  \label{fact:nmc}
  Let $T$ be the theory of rich structures and assume that every $\omega$-saturated model of $T$ is rich. Suppose that for all $\X\vDash T$ and $\bc\se X$, there exists an existential $L$-formula $\tau_{\bc}(\bz)$ such that
  \begin{enumerate}
    \item $\X\vDash\tau_{\bc}(\bc)$, and\label{it:nmc1}
    \item for all $\Y\vDash T$ and $\bc'\se Y$, if $\Y\vDash\tau_{\bc}(\bc')$ then $\delta(\bc')\leq\delta(\bc)$.\label{it:nmc2}
  \end{enumerate}
  Then $T$ is near model complete.
\end{fact}

In our case, the predimension function $\delta$ is defined for pairs, and here is a verification of Fact \ref{fact:nmc} for $T=T_0^{\rich}$: Let $L^*$ be the expansion of $L$ by predicates for existentially definable sets, and let $T^*$ be the corresponding theory extending $T$. We prove quantifier elimination of $T^*$ by embedding test. Let $\X,\Y$ be $\omega$-saturated models of $T$, and let $\ba_1\se X,\ba_2\se Y$ such that $\ba_1\to\ba_2$ is a partial $L^*$-isomorphism. We want to extend it to a partial isomorphism in $\F(\X,\Y)$. Take $A_1:=[\ba_1]^\X$ which is strong in $X^2$. The quantifier free type of $\ag{A_1}$ over $\ba_1$ is witnessed by $T^*$, and hence by saturation, there exists $A_2\se Y^2$ such that $\ag{A_1}\to\ag{A_2}$ is a partial isomorphism extending $A_1\to A_2$ and $\ba_1\to\ba_2$. It suffices to show that $A_2$ is strong in $Y^2$. Let $\bb_1$ be a $\cl$-basis of $A_1$, and let $\bb_2$ be the image of $\bb_1$, a $\cl$-basis of $A_2$. Suppose on the contrary that there is $\bc\se Y^2$ such that $\delta(\bc\bb_2)<\delta(\bb_2)$. Then $\X\vDash\exists\bz,\tau_{\bc\bb_2}(\bz\bb_1)$ and $\Y\vDash\tau_{\bb_1}(\bb_2)$. Thus, there exists $\bc'\se X^2$ such that
\[\delta(\bc'\bb_1)\leq\delta(\bc\bb_2)<\delta(\bb_2)\leq\delta(\bb_1),\]
contradicting the fact that $\cl(\bb_1)$ is strong in $X^2$.

To find such formulas, we first work on multiplicatively independent tuples.

\begin{lemma}
  \label{le:mi}
  Let $\X$ be a model of $T_0^{\rich}$, let $c_1,...,c_n$ be multiplicatively independent in $G^\X$, and let $V\se R^{2n}$ be the algebraic locus of $(c_1,...,c_n)$ over $\bQ$. Define $\tau_{\bc}(\bz)=\tau_{\bc}(x_1,y_1,...,x_n,y_n)$ as the conjunction of the following formulas:
  \begin{enumerate}
    \item $(x_1,y_1,...,x_n,y_n)\in V\cap G^n$,
    \item $\dim^K(\widetilde{V}(\bx+i\by,-))=r_V$, \label{it:tau2}
    \item $\bigwedge_{B\in\mu(\widehat{V})}(x_1,y_1,...,x_n,y_n)\notin B$. \label{it:tau3}
  \end{enumerate}
  Then $\tau_{\bc}$ meets the requirement of Fact \ref{fact:nmc}.
\end{lemma}
\begin{proof}
  The idea is similar to the proof of Lemma \ref{le:c0}. First, observe that $\bc$ satisfies $\tau_{\bc}$. Let us show that statement \eqref{it:nmc2} in Fact \ref{fact:nmc} holds. $\delta(\bc)$ can be expressed by:
  \[\delta(\bc)=\trd(\bc)-\md(\bc)=\dim(V)-n.\]
  
  Suppose $\Y\vDash T_0^{\rich}$ and $\bb\se Y^2$ such that $\Y\vDash\tau_{\bc}(\bb)$. By \eqref{it:tau2} of $\tau_{\bc}$, we have
  \begin{equation}
    \label{eq:trdbd}
    \trd(\bb/\hb)\leq\dim^K(\widetilde{V}(\hb,-))=\dim(V)-\dim^K(\widehat{V}).
  \end{equation}
  \begin{claim}
    $\delta(\bb)\leq\dim(V)-n.$
  \end{claim}
  \begin{pfcl}
    Suppose on the contrary that the following holds.
    \begin{equation}
    \label{eq:contrary}
      \trd(\bb)-\md(\bb)>\dim(V)-n
    \end{equation}
    Then the following inequality shows that $\hb$ is multiplicatively dependent:
    \[\md(\hb)=\md(\bb)<\trd(\bb)-\dim(V)+n\leq n.\]

    Let $B$ be an algebraic subgroup containing $\hb$ with $\dim^K(B)=\md(\hb)$. Let $S$ be an irreducible component of $\widehat{V}\cap B$ containing $\hb$. Combining \eqref{eq:trdbd} and \eqref{eq:contrary}, we see that $S$ is atypical:
    \[\dim^K S\geq\trd\hb\geq\dim^K\widehat{V}+\trd\bb-\dim V>\dim^K\widehat{V}+\dim^K B-n.\]
    By Fact \ref{fact:WCIT} and Lemma \ref{le:id}, there is $B\in\mu(\widehat{V})$ such that $\bb\in B$, contradicting \eqref{it:tau3} of $\tau_{\bc}$.\hfill$\blacksquare$
  \end{pfcl}

  Therefore, $\delta(\bb)\leq\delta(\bc)$. We conclude that $\tau_{\bc}$ satisfies the desired property.
\end{proof}

Let us now construct the formula for a general tuple.
\begin{lemma}
  \label{le:gn}
  Let $\X\vDash T_0^{\rich},\,\bc\se X^2$, let $\bb$ be a $\cl$-basis of $\cl(\bc)\cap G^\X$, and let $W(\by,\bz)\se(R^2)^{|\bb|+|\bc|}$ be the algebraic locus of $(\bb,\bc)$ over $\bQ$. Let $\tau_{\bb}$ be defined as in Lemma \ref{le:mi}. Define $\tau_{\bc}(\bz)$ to be the following formula:
  \[\exists\,\by,\;\;\tau_{\bb}(\by)\,\bigwedge\,W(\by,\bz)\,\bigwedge\,\dim(W(\by,-))=\trd(\bc/\bb).\]
  Then $\tau_{\bc}$ meets the requirement of Fact \ref{fact:nmc}.
\end{lemma}
\begin{proof}
  Statement \eqref{it:nmc1}, i.e., $\X\vDash\tau_{\bc}(\bc)$, is obvious.

  Let $\Y\vDash T_0^{\rich},\,\bc'\se Y^2$ satisfying $\Y\vDash\tau_{\bc}(\bc')$. We need to show that $\delta(\bc')\leq\delta(\bc)$.

  Let $\bb'\se Y^2$ be an instance of $\by$ in $\tau_{\bc}(\bc')$. Lemma \ref{le:mi} implies $\delta(\bb')\leq\delta(\bb)$. Notice that $W(\bb',\bc')$ gives $\bb'\se\cl(\bc')$. Moreover,
  \[\delta(\bc'/\bb')\leq\trd(\bc'/\bb')\leq\dim(W(\bb',-))=\trd(\bc/\bb)=\delta(\bc/\bb).\]
  Thus,
  \[\delta(\bc')=\delta(\bc'/\bb')+\delta(\bb')\leq\delta(\bc/\bb)+\delta(\bb)=\delta(\bc).\qedhere\]
\end{proof}

By Fact \ref{fact:nmc}, we arrive at the conclusion.

\begin{proposition}
  $T_0^{\rich}$ is near model complete.\qed
\end{proposition}

\subsection{O-minimal open core}

We want to study open sets definable in an expansion. Let $\M$ be an expansion of a dense linear order without endpoints. The open core of $\M$ is a relational structure with underlying set $M$, and for each open definable set $U\se M^n$ a predicate belonging to its language. Let $\M^*$ be an expansion of $\M$ by language. Boxall and Hieronymi \cite{Bo12} found a criterion for the expansion to preserve its open core. Assume that $\M^*$ is sufficiently saturated and strongly homogeneous.

\begin{fact}[{\cite[Corollary 3.1]{Bo12}}]
  \label{fact:oc}
  Let $C$ be a small set in $M$. Suppose that for every $n\in\bN$ there is a set $D_n\se M^n$ such that the following conditions hold:
  \begin{enumerate}
    \item $D_n$ is dense in $M^n$,\label{it:open1}
    \item for every $\bx\in D_n$ and every open set $U\in M^n$, if $\tp_{\M}(\bx|C)$ is realized in $U$ then it is realized in $U\cap D_n$,\label{it:open2}
    \item for every $\bx\in D_n$, the conjunction of $\tp_{\M}(\bx|C)$ and $\bx\in D_n$ implies $\tp_{\M^*}(\bx|C)$.\label{it:open3}
  \end{enumerate}
  Then every open set definable over $C$ in $\M^*$ is definable over $C$ in $\M$.
\end{fact}

In our context, $\M^*=(\M,G)$ is a rich structure in $\C_0$, and $\M$ is a real closed field. Without loss of generality, assume that $C$ is a finite tuple $\bc$. By Fact \ref{fact:min}, extend $\bc$ such that $(\ag{\bc},G^\ag{\bc})\leq\M^*$. Let $D_n$ be the following set:
\[D_n=\{\,\bx\in M^n\mid\forall\ba\se G,\,\trd(\ba\bx/\bc)-\md(\ba/G^\ag{\bc})\geq n\,\}.\]
In other words, $D_n$ consists of all $n$-tuples $\Cl$-independent over $\ag{\bc}$.

We want to verify the three conditions in Fact \ref{fact:oc}. If every open rectangle in $M^n$ intersects with $D_n$, then \eqref{it:open1} and \eqref{it:open2} are immediate.

\begin{lemma}
  Let $U$ be an open rectangle in $M^n$. Then $U\cap D_n\neq\emptyset$.
\end{lemma}
\begin{proof}
  Let $u_1,v_1,...,u_n,v_n\in M$ be such that $U=(u_1,v_1)\times...\times(u_n,v_n)$. Let $\bc'$ be a tuple containing $\bc$ and all $u_i,v_i$'s such that $(\ag{\bc'},G^\ag{\bc'})\leq\M^*$. There exist a structure $\A\in\Cf_0$ and $a_1,...,a_n\in A$ such that:
  \begin{itemize}
    \item $A=\ag{\bc',a_1,...,a_n}$,
    \item $\trd(a_1,...,a_n/\bc')=n$,
    \item $G^\A=G^\ag{\bc'}$, and
    \item $a_i\in(u_i,v_i)$ for all $i=1,...,n$.
  \end{itemize}
  Since $\M^*$ is rich, there is a strong embedding $\iota:\A\to\M^*$ extending $\ag{\bc'}\leq\M^*$. Then $\iota(a_1),...,\iota(a_n)$ are $\Cl$-independent over $\ag{\bc'}$, and therefore,
  \[(\iota(a_1),...,\iota(a_n))\in D_n.\qedhere\]
\end{proof}

We show condition \eqref{it:open3} using the back-and-forth system in Subsection \ref{sec:bf}.

\begin{lemma}
  Let $\bx,\by\in D_n$ be such that $\by$ satisfies $\tp_{\M}(\bx|\bc)$. Then there is a partial isomorphism $f\in\F(\M^*,\M^*)$, i.e.,
  \[f:\X_1\xrightarrow{\sim}\X_2,\text{ where }\X_i\leq\M^*,\X_i\in\Cf_0,\]
  fixing $\bc$ and sending $\bx$ to $\by$.
\end{lemma}
\begin{proof}
  We observe that, because $\bx,\by\in D_n$,
  \[G^\ag{\bc}=G^\ag{\bc\bx}=G^\ag{\bc\by}.\]
  Therefore, the map $f:\ag{\bc\bx}\to\ag{\bc\by}$ fixing $\bc$ and sending $\bx$ to $\by$ is a partial isomorphism. Let $\X_1$ and $\X_2$ be the substructures with underlying sets $\ag{\bc\bx}$ and $\ag{\bc\by}$ respectively. Check that $\X_i\leq\M^*,\X_i\in\Cf_0$, and hence $f\in\F(\M^*,\M^*)$.
\end{proof}

Partial isomorphisms in a back-and-forth system are elementary. Thus, condition \eqref{it:open3} holds. Applying Fact \ref{fact:oc}, we conclude that $\M^*$ defines no new open set, and hence has an o-minimal open core.

\begin{corollary}
  Let $\U$ be a model of $T_0^{\rich}$. Then every open set definable in $\U$ is semialgebraic.
\end{corollary}

\subsection{Strong dependence}

Let us recall some definitions and facts about dp-rank. For more details, we refer to Simon's book \cite[Chapter 4]{Si15}.

\begin{definition}
  Let $p$ be a partial type over a set $A$, and let $\kappa$ be a cardinal. We say $\dprk(p,A)<\kappa$ if for every family $(I_t\mid t<\kappa)$ of mutually indiscernible sequences over $A$ and $\bb\vDash p$, there is $t<\kappa$ such that $I_t$ is indiscernible over $A\bb$.
\end{definition}

For a tuple $\ba$, we simply write $\dprk(\ba/A)$ for $\dprk(\tp(\ba|A),A)$.

\begin{fact}
  $\dprk(p,A)<\kappa$ if and only if the following holds:

  For every family $(I_t\mid t<\kappa)$ of mutually indiscernible sequences over $A$ and $\bb\vDash p$, there is $t<\kappa$ such that all members of $I_t$ have the same type over $A\bb$.
\end{fact}

\begin{fact}
  A theory $T$ is NIP if and only if for every type $p$ and set $A$ there is some $\kappa$ such that $\dprk(p,A)<\kappa$.
\end{fact}

\begin{definition}
  An NIP theory $T$ is \emph{strongly dependent} if $\dprk(x=x,\emptyset)<\aleph_0$.
\end{definition}

\begin{fact}\label{fa:dpad}
  Let $\U$ be a monster model and $\ba,\bb\se\U$. Let $A$ be a small set and $\kappa_1,\kappa_2$ be two cardinals such that $\dprk(\bb/A)<\kappa_1$ and $\dprk(\ba/A\bb)<\kappa_2$. Then $\dprk(\ba\bb/A)<\kappa_1+\kappa_2-1$.
\end{fact}

We want to compute the dp-rank in rich structures. Fix a sufficiently saturated and strongly homogeneous model $\U\vDash T^{\rich}_0$. Our first observation is that the definable closure is rcl-generated by the hull.

\begin{proposition}\label{prop:dcl}
  The definable closure in $\U$ is given by $\dcl(A)=\ag{[A]}$ for every small $A\se U$.
\end{proposition}
\begin{proof}
  It suffices to show $\dcl(\bx)=\ag{[\bx]}$ for every $\bx\se U$. Let $\X$ denote the substructure with underlying set $\ag{[\bx]}$.

  $(\se):$ Take $a\in U\setminus X$. Then there is some $\Y\leq\U$ with $a\in Y$ and $\X\leq\Y$. Lemma \ref{le:inf} implies that $a$ cannot be definable.

  $(\supseteq):$ By Fact \ref{fact:min}, there is $\by\se G^\U$ such that $[\bx]=\cl(\by)$. It suffices to show that $\by\in\dcl(\bx)$. Every automorphism fixing $\bx$ fixes $[\bx]$, and hence $\by$ has at most countably many conjugates over $\bx$. If $\by\notin\operatorname{acl}(\bx)$ and $\U$ is $\kappa$-saturated, then $\by$ has at least $\kappa$ many conjugates. Therefore, we must have $\by\in\operatorname{acl}(\bx)=\dcl(\bx)$.
\end{proof}

\begin{lemma}\label{le:dp}
  Let $\X,\Y\leq\U$ be structures in $\Cf_0$ such that $\X\leq\Y$ is minimal prealgebraic and $\trd(Y/X)=n$. Let $((\Z_i^j)_{i<\omega})_{j<n+1}$ be a family of mutually indiscernible sequences over $\X$ such that $\X\leq\Z_i^j\leq\U$ and $\Z_i^j\in\Cf_0$ for each $i<\omega,j<n+1$. Then there is $j<n+1$ such that $(\Z_i^j)_{i<\omega}$ have the same type over $\Y$.
\end{lemma}
\begin{proof}
  Compare $Y\cap Z_i^j$ and $X$. There are two cases:
  \begin{enumerate}[leftmargin=*]
    \item $Y\cap Z_t^u\supsetneqq X$ for some $t<\omega,u<n+1$.
    
    Let $\Y'$ denote the substructure with underlying set $Y\cap Z_t^u$. Fact \ref{fact:str} implies $\Y'\leq\U$. Since $\X\leq\Y$ is minimal, we have $Y\se Y'$. By definition, an indiscernible sequence $(\Z_i^j)_{i<\omega}$ with $j\neq u$ is indiscernible over $\Z_t^u$ and hence over $\Y$.

    \item $Y\cap Z_i^j=X$ for all $i,j$.
    
    Lemma \ref{le:inters} implies $Y\forkindep^{\rcl}_X Z_i^j$ for all $i,j$. In real closed fields, dp-rank equals transcendence degree. Therefore, there is $j<n+1$ such that the types in the reduct $\tp_{\rcf}(Z_i^j|Y)$ are the same for all $i<\omega$. Let $\Z_i'$ denote the substructure with underlying set $\ag{Y\cup Z_i^j}$. Since $\X\leq\Y$ is prealgebraic, we have
    \[\md(G^{\Z_i'}/G^{\Z_i^j})\leq\trd(Z_i'/Z_i^j)=\trd(Y/X)=\md(G^\Y/G^\X)=\md(G^\Y/G^{\Z_i^j}).\]
    Therefore, $G^{\Z_i'}=G^\Y\cdot G^{\Z_i^j}$ and $\Z_i^j\leq\Z_i'$ prealgebraic, resulting in $\Z_i'\leq\U$. Thus, we obtain partial isomorphisms $\Z_{i_1}'\to\Z_{i_2}'$ for all $i_1,i_2<\omega$ in the back-and-forth system. In conclusion, types $\tp(Z_i^j|Y)$ are the same for all $i<\omega$.\qedhere
  \end{enumerate}
\end{proof}

\begin{proposition}
  Let $\X\leq\U$ be a structure in $\Cf_0$. Then $\dprk(b/X)<\aleph_0$ for every $b\in U$.
\end{proposition}
\begin{proof}
  Let $((I_i^j)_{i<\omega})_{j<\omega}$ be a family of mutually indiscernible sequences over $X$. By Lemma \ref{le:hull0}, let $\Z_i^j$ be the strong substructure with underlying set $\ag{[XI_i^j]}$. Proposition \ref{prop:dcl} implies $\dcl(XI_i^j)=\ag{[XI_i^j]}$. Thus, $((\Z_i^j)_{i<\omega})_{j<\omega}$ is also mutually indiscernible over $\X$. Again by Lemma \ref{le:hull0}, let $\Y$ be the strong substructure with underlying set $\ag{[Xb]}$. There are two cases for the type of $b$ and the type of $\X\leq\Y$:
  \begin{enumerate}[leftmargin=*]
    \item $b\in\Cl(X)$ and $\X\leq\Y$ prealgebraic.
    
    Let $\X=\Y_0\leq\Y_1\leq...\leq\Y_n=\Y$ be a chain of minimal prealgebraic extensions. Lemma \ref{le:dp} implies $\dprk(Y_{i+1}/Y_i)<\trd(Y_{i+1}/Y_i)+1$. Combining this with Fact \ref{fa:dpad}, we get $\dprk(Y/X)<\trd(Y/X)+1$.

    \item $b\notin\Cl(X)$ and $\X\leq\Y$ purely transcendental.
    
    In this case, $Y=\ag{Xb}$ and $G^\Y=G^\X$. If $Z_i^j\cup\{b\}$ is strong in $U$ for all $i,j$, then by the dp-rank of RCF, there is $j<\omega$ such that $\tp_{\rcf}(Z_i^j|Y)$ are the same for all $i<\omega$, and the back-and-forth system implies $\tp(Z_i^j|Y)$ are the same for all $i<\omega$.

    Suppose there is some $t,u$ such that $Z_t^u\cup\{b\}$ is not strong in $U$. Then $\Y':=\ag{[Z_t^ub]}$ must be prealgebraic over $\Z_t^u$. Now $((\Z_i^j)_{i<\omega})_{j\neq u}$ is mutually indiscernible over $\Z_t^u$. Apply Lemma \ref{le:dp} substituting $\X,\Y$ with $\Z_t^u,\Y'$ to acquire some $j\neq u$ such that $(\Z_i^j)_{i<\omega}$ have the same type over $Y'$.\qedhere
  \end{enumerate}
\end{proof}

\begin{corollary}
  $T_0^{\rich}$ is strongly dependent.
\end{corollary}

\section{Model}
\label{ch:md}

In this section, we present a concrete model of the theory of rich structures. Let $(\overline{\bR},G)$ be the expansion of the real ordered field by a subset
\[G=\exp(\epsilon\bR+Q)=\{e^{\epsilon t+s}\in\bR^2\mid t\in\bR,s\in Q\},\]
where $Q\se\bR$ is a non-trivial finite dimensional $\bQ$-vector space, and $\epsilon=1+\beta i$ for some $\beta\in\bR^\times$. To show Theorem \ref{thm:md}, we verify the two components of $T_\B^{\rich}$, which are $T_\B$ and EC.

\subsection{Predimension inequality}
We want to prove that $(\overline{\bR},G)$ is a model of $T_\B$, for a suitable $\B$. We adopt a similar approach to Caycedo's \cite[Section 6.2]{Ca11}.

\begin{lemma}
  \label{le:lb}
  Assume SC$_K$ (Conjecture \ref{con:SCK}) holds for $K=\bQ(\beta i)$. Then for all $\ba=\bigl((p_1,q_1),...,(p_n,q_n)\bigr)\se(\bR^2)^\times$, we have
  \[\delta(\ba)=\trd(\ba)-\md(\cl(\ba)\cap G)\geq-3\ld_\bQ(Q)-\trd(K).\]
\end{lemma}
\begin{proof}
  We may assume that $\ba\se G$ and is multiplicatively independent. Let $\by=\bp+i\bq=(p_1+iq_1,...,p_n+iq_n)\in(\bC^\times)^n$. Take $\bx\in\bC^n$ such that $e^{\bx}=\by$. The following equalities are immediate:
  \begin{gather*}
    \trd(\ba)=\trd(\bp\bq)=\trd(\by\by^c),\\
    n=\md(\cl(\ba)\cap G)=\md(\by)=\ld_\bQ(\bx).
  \end{gather*}
  Applying SC$_K$ to $\bx\bx^c$, we have
  \[\trd(\by\by^c)\geq\ld_\bQ(\bx\bx^c)-\ld_K(\bx\bx^c)-\trd(K).\]
  Therefore, it suffices to prove that $\ld_\bQ(\bx\bx^c)-\ld_K(\bx\bx^c)\geq n-3\ld_\bQ(Q)$.

  Caycedo has proved the desired bound. The following two inequalities are from the proof of \cite[Lemma 6.2.1]{Ca11}. Write $\bx=\epsilon\bt+\bs$, where $\bt\se\bR,\bs\se Q$.
  \begin{claim}
    $\ld_K(\bx\bx^c)\leq n+\ld_\bQ(Q)$.
  \end{claim}
  \begin{pfcl}
    Observe that $\epsilon$ and $\epsilon^c$ are $K$-linearly dependent. Indeed,
    \[\frac{\epsilon}{\epsilon^c}=\frac{1+\beta i}{1-\beta i}\in\bQ(\beta i)=K.\]
    As a result,
    \[\ld_K(\bx\bx^c)\leq\ld_K(\epsilon\bt,\bs)\leq n+\ld_\bQ(Q).\tag*{$\blacksquare$}\]
  \end{pfcl}
  \begin{claim}
    $\ld_\bQ(\bx\bx^c)\geq 2n-2\ld_\bQ(Q)$.
  \end{claim}
  \begin{pfcl}
    Denote $r=\ld_\bQ(Q)$. Without loss of generality, we can assume $x_{i}=\epsilon t_i$ for all $i>r$. Indeed, we may apply $\bQ$-linear transformations to $\bx\bx^c$ to reduce to this case. Recall that $(x_{r+1},...,x_n)$, and hence $(t_{r+1},...,t_n)$, are $\bQ$-linearly independent. Thus,
    \[\ld_\bQ(\bx\bx^c)\geq\ld_\bQ(\epsilon t_{r+1},...,\epsilon t_n,\epsilon^c t_{r+1},...,\epsilon^c t_n)=2(n-r).\tag*{$\blacksquare$}\]
  \end{pfcl}
  In conclusion,
  \[\delta(\ba)=\trd(\by\by^c)-n\geq-3\ld_\bQ(Q)-\trd(K).\qedhere\]
\end{proof}

\begin{corollary}
  Assume SC$_K$ holds for $K=\bQ(\beta i)$. Then there is $\bb\se G$ such that $\bb$ is multiplicatively independent and $\delta(\ba/\bb)\geq 0$ for all $\ba\se(\bR^2)^\times$.
\end{corollary}
\begin{proof}
  Since Lemma \ref{le:lb} provides a lower bound of $\delta$, by Fact \ref{fact:min}, there is a finite tuple $\bb$ achieving the minimum of $\delta$. By removing redundant elements in $\bb$, we may assume that $\bb\se G$ and $\bb$ is multiplicatively independent.
\end{proof}

With $\bb$ in the above corollary, let $\B\se(\overline{\bR},G)$ with underlying set $\ag{\bQ,\bb}$. We observe that $(\overline{\bR},G)\vDash T_\B$.

\subsection{EC-property}
This subsection is devoted to the proof of the EC-property for $(\overline{\bR},G)$, i.e., every rotund block intersects with $G^n$. Our argument is inspired by Caycedo and Zilber's \cite{Ca14}.

Let $V\se\bR^{2n}$ be a rotund block defined over $\bc\se\bR$. Denote
\[X=\{(\bt,\bs)\in\bR^{2n}\mid\exp(\epsilon\bt+\bs)\in V\}.\]
It suffices to show that $X\cap(\bR^n\times Q^n)$ is non-empty. Let $\pi:\bR^{2n}\to\bR^n$ be the projection onto the last $n$ coordinates. We prove that $\pi(X)$ contains an open subset of $\bR^n$.

Let $\R$ denote the expansion of $(\bR,<,+,-,\cdot,0,1)$ by
\begin{itemize}
  \item constants $\bc$, and
  \item functions $x\mapsto e^x$ and $x\mapsto\sin\beta x$ defined on bounded intervals with rational endpoints.
\end{itemize}
As a reduct of $\bR_{an}$, the structure $\R$ is o-minimal. We write $\dim$ for the topological dimension in $\R$. Define the following function.
\begin{align*}
  f:\;\;\;\bR^n\times\bR^n\;&\longrightarrow\;\;\;\bR^{2n}\\
  (\bt,\,\bs)\;\;\;&\longmapsto\;\exp(\epsilon\bt+\bs)
\end{align*}
Clearly, $f|_U$ is definable in $\R$ for every open rectangle $U\se\bR^{2n}$ with rational endpoints.
\begin{lemma}
  There is an open subset $U\se\bR^{2n}$ such that
  \begin{itemize}
    \item $U$ and $f|_U$ are definable in $\R$,
    \item $f|_U$ is a homeomorphism onto its image,
    \item $f(U)\cap V$ is still a rotund block.
  \end{itemize}
\end{lemma}
\begin{proof}
  Let $\ba$ be a generic point of $V$. Since $V$ is rotund, we have $\ba\se(\bR^2)^\times$ and $f^{-1}(\ba)$ is discrete. There is an open rectangle $E\se\bR^{2n}$ with rational endpoints such that $\ba\in E$ and $f$ restricts to a homeomorphism $U\xrightarrow{\sim}E$ for some $U$ bounded open in $\bR^{2n}$. By shrinking $E$, we may assume $E\cap V$ to be connected. Observe that $U$ is definable in $\R$. Indeed, we may define $f$ on an open rectangle containing $U$, and then take the preimage of $E$. Since $\ba$ is generic in $E\cap V$, Lemma \ref{le:dimNash} implies $V^{\zar}=(E\cap V)^{\zar}$ and $E\cap V$ is rotund.
\end{proof}
Fix such an open subset $U$. From now on, let us simply write $f$ for $f|_U$. Define $V'=f(U)\cap V$ and $X'=U\cap X$. We see that $f$ restricts to a homeomorphism $X'\xrightarrow{\sim}V'$. Consequently, $\dim(X')=\dim(V')=n$. Our key observation is the following.
\begin{lemma}
  \label{le:fb}
  For every generic point $(\bt^*,\bs^*)\in X'$, the fiber $X'_{\bs^*}$ has dimension $0$.
\end{lemma}
We will go back to the proof of this lemma later. Let us now continue with the proof of the EC-property.
\begin{lemma}
  \label{le:pi}
  The set $\pi(X')$ has dimension $n$.
\end{lemma}
\begin{proof}
  The following set is definable in $\R$:
  \[X'_0=\{\,(\bt,\bs)\in X'\mid\dim X'_{\bs}=0\,\}.\]
  By Lemma \ref{le:fb}, every generic point of $X'$ is contained in $X'_0$, or in other words,
  \[\dim(X'\setminus X'_0)<\dim(X').\]
  Therefore, $\dim(X'_0)=\dim(X')=n$. The dimension of fibers in an o-minimal structure satisfies:
  \[\dim\pi(X'_0)=\dim(X'_0)-0=n.\qedhere\]
\end{proof}
\begin{proposition}
  The set $V\cap G^n$ is non-empty. Hence, $(\overline{\bR},G)$ has the EC-property.
\end{proposition}
\begin{proof}
  Lemma \ref{le:pi} implies that $\pi(X')$ contains an open subset of $\bR^n$. Then $Q^n$ must intersect with $\pi(X')$, because $Q^n$ is dense in $\bR^n$. Take an element $(\bt,\bs)\in X'\cap(\bR^n\times Q^n)$. Then $\exp(\epsilon\bt+\bs)\in V\cap G^n$.
\end{proof}

\subsubsection{Unlikely intersection}

Our goal is Lemma \ref{le:fb}.

Take any generic point $(\bt^*,\bs^*)\in X'$, and hence $\exp(\epsilon\bt^*+\bs^*)$ is generic in $V'$. We want to show that $X'_{\bs^*}$ has dimension $0$. Let $B:=U_{\bs^*}$. It is an open neighborhood of $\bt^*$ in $\bR^n$, and
\[f(X'_{\bs^*}\times\{\bs^*\})=\exp(\epsilon B+\bs^*)\cap V'.\]
It suffices to prove the following.

\begin{lemma}
  The set $\exp(\epsilon B+\bs^*)\cap V'$ has dimension $0$.
\end{lemma}
\begin{proof}
  Suppose on the contrary that
  \[\dim X'_{\bs^*}=\dim\left(\exp(\epsilon B+\bs^*)\cap V'\right)>0.\]
  Since $\bt^*$ is generic, we have that $X'_{\bs^*}$ contains an analytic curve through $\bt^*$, which can be parametrized by $\bt(x):(-1,1)\to\bR^n$ with $\bt(0)=\bt^*$. Consequently,
  \[\exp(\epsilon\bt(x)+\bs^*)\in V'\text{ for all }x\in(-1,1).\]
  Since $\bt(x)$ is real-analytic and non-constant, it can be expressed by $(t_1,...,t_n)\se\bR[[x]]$ with $\rk J(t_1,...,t_n)>0$.

  In order to apply Fact \ref{fact:ax}, we need a $\bQ$-linearly independent set modulo $\bC$. Define $t_i'(x)=t_i(x)-t_i(0)$, i.e, $t_i(x)=t_i'(x)+t_i^*$. Take a $\bQ$-linear basis of $\{t_1',...,t_n'\}$, say $\{t_{\sigma(1)}',...,t_{\sigma(k)}'\}$ for some permutation $\sigma$ and $1\leq k\leq n$. To simplify our notation, let us assume that the first $k$ terms of $t_1',...,t_n'$ provide a basis. Since $t_1',...,t_n'\in x\bR[[x]]$, we see that $t_1',...,t_k',i\beta t_1',...,i\beta t_k'$ are $\bQ$-linearly independent modulo $\bC$. By Fact \ref{fact:ax},
  \begin{equation}
    \label{eq:ax}
    \trd_\bC(t_1',...,t_k',i\beta t_1',...,i\beta t_k',e^{t_1'},...,e^{t_k'},e^{i\beta t_1'},...,e^{i\beta t_k'})>2k.
  \end{equation}
  On the other hand, we also have the following relations:
  \[\trd_\bC(\bt',i\beta\bt')=\trd_\bC(\bt')\leq k,\]
  \[\trd_\bC(e^{\bt'},e^{i\beta\bt'})=\trd_\bC\left(\operatorname{Re}(e^{\epsilon\bt+\bs^*}),\operatorname{Im}(e^{\epsilon\bt+\bs^*})\right),\]
  where $\operatorname{Re}(-)$ and $\operatorname{Im}(-)$ denote the real and imaginary parts.
  \begin{claim}
    $\trd_\bC\left(\operatorname{Re}(e^{\epsilon\bt+\bs^*}),\operatorname{Im}(e^{\epsilon\bt+\bs^*})\right)\leq k$.
  \end{claim}
  \begin{pfcl}
    There is a matrix $M\in\bZ^{(n-k)\times n}$ of full rank such that $\bt'\in\ker(M)$. Indeed, for each $j\leq n-k$, we have that $t_{k+j}'=a_{1j} t_1'+...+a_{kj} t_k'$ for some $a_{ij}\in\bQ$. Let $M$ be a multiple of the following matrix:
    \[\begin{bmatrix}
      a_{11}&...&a_{k1}&-1&0&...&0\\
      a_{12}&...&a_{k2}&0&-1&...&0\\
      ...&...&...&...&...&...&...\\
      a_{1(n-k)}&...&a_{k(n-k)}&0&0&...&-1
    \end{bmatrix}\]
    Clearly, $M\bt'=0$ and $\rk(M)=n-k$. Let us write $(-)^M$ for the map induced by $M$ with respect to complex power, and write $A$ for the algebraic subgroup defined by $\bz^M=\mathds{1}$. Because $V'$ is rotund, $\dim(V'^M)\geq n-k$. Since $e^{\epsilon\bt^*+\bs^*}$ is a generic point of $V'$, the dimension of its fiber satisfies:
    \[\dim\left(V'\cap(e^{\epsilon\bt^*+\bs^*}A)\right)=\dim(V')-\dim(V'^M)\leq n-(n-k)=k.\]
    Notice that $e^{\epsilon\bt(x)+\bs^*}$ is contained in $V'\cap(e^{\epsilon\bt^*+\bs^*}A)$ for all $x\in(-1,1)$. Hence,
    \[\trd_\bC\left(\operatorname{Re}(e^{\epsilon\bt+\bs^*}),\operatorname{Im}(e^{\epsilon\bt+\bs^*})\right)\leq\dim\left(V'\cap(e^{\epsilon\bt^*+\bs^*}A)\right)\leq k.\tag*{$\blacksquare$}\]
  \end{pfcl}

  In conclusion, we have the following inequality, contradicting \eqref{eq:ax}:
  \[\trd_\bC(\bt',i\beta\bt',e^{\bt'},e^{i\beta\bt'})\leq\trd_\bC(\bt')+\trd_\bC\left(\operatorname{Re}(e^{\epsilon\bt+\bs^*}),\operatorname{Im}(e^{\epsilon\bt+\bs^*})\right)\leq 2k.\qedhere\]
\end{proof}

This finishes the proof of Lemma \ref{le:fb}.

\bibliography{paper}
\bibliographystyle{plain}

\end{document}